\documentclass[preprint,number,3p]{elsarticle}

\usepackage{graphicx}

\usepackage{amssymb,amsthm,amsmath}

\usepackage{subfig}
\usepackage{float}
\usepackage{multirow}
\usepackage{longtable,lscape}
\usepackage{cases}
\usepackage{bm,diagbox}
\usepackage[colorlinks=true]{hyperref}
\usepackage{url}
\usepackage{enumitem}

\usepackage[usenames, dvipsnames]{color}

\usepackage{fancyhdr}
\newtheorem{theorem}{Theorem}[section]

\theoremstyle{definition}
\newtheorem{definition}[theorem]{Definition}
\newtheorem{remark}{Remark}

\usepackage{titlesec}
\titleformat{\paragraph}
{\normalfont\normalsize\bfseries}{\theparagraph}{1em}{}
\titlespacing*{\paragraph}
{0pt}{3.25ex plus 1ex minus .2ex}{1.5ex plus .2ex}

\numberwithin{equation}{section}

\journal{DCDS-B}

\usepackage{algorithm}
\usepackage{algorithmicx}
\usepackage{algpseudocode}
\algrenewcommand\algorithmicrequire{\textbf{Input:}}
\algrenewcommand\algorithmicensure{\textbf{Output:}}

\usepackage{caption}
\newcommand{\norm}[1]{\left\lVert#1\right\rVert}
\DeclareMathOperator*{\argmin}{argmin}

\begin{document}

\begin{frontmatter}



\title{{\vspace*{-4cm} \small \color{red}Please cite this paper as: Huang, Kuang, et al. ``A game-theoretic framework for autonomous vehicles velocity control: Bridging microscopic differential games and macroscopic mean field games.'' Discrete \& Continuous Dynamical Systems-B 22.11 (2017): 0.  \vspace*{4cm}} \\ A Game-Theoretic Framework for Autonomous Vehicles Velocity Control: Bridging Microscopic Differential Games and Macroscopic Mean Field Games}


\author[cuAPAM]{Kuang Huang}

\author[cu,dsi]{Xuan Di\corref{cor}}
\ead{sharon.di@columbia.edu}

\author[cuAPAM,dsi]{Qiang Du}
\author[cuCS]{Xi Chen}

\cortext[cor]{Corresponding author. Tel.: +1 212 853 0435;}

\address[cuAPAM]{Department of Applied Physics and Applied Mathematics, Columbia University}
\address[cu]{Department of Civil Engineering and Engineering Mechanics, Columbia University}
\address[dsi]{Data Science Institute, Columbia University}
\address[cuCS]{Department of Computer Science, Columbia University}

\begin{abstract}
  This paper 
  proposes an efficient computational framework for longitudinal velocity control of a large number of autonomous vehicles (AVs) 
  and develops a traffic flow theory for AVs.
  Instead of hypothesizing explicitly how AVs drive,
  our goal is
  to design future AVs as rational, utility-optimizing agents
  that continuously select optimal velocity over a period of planning horizon. 
  With a large number of interacting AVs,
  this design problem can become computationally intractable.
  This paper aims to tackle such a challenge by employing mean field approximation and deriving a
  mean field game (MFG) as the limiting differential game with an infinite number of agents.
  The proposed micro-macro model
  allows one to define individuals on a microscopic level as utility-optimizing agents while translating rich microscopic behaviors to macroscopic models.
  Different from existing  studies on the application of MFG to traffic flow models,  
  the present study
  offers a systematic framework to apply MFG to autonomous vehicle velocity control. 
  The MFG-based AV controller is shown to mitigate traffic jam faster than the LWR-based controller.
  MFG also embodies classical traffic flow models with behavioral interpretation,
  thereby providing a new traffic flow theory for AVs.
\end{abstract}


\end{frontmatter}

\section{Introduction}\label{sec:motive}

\subsection{Problem statement}

When all human-driven vehicles (HVs) on public roads are replaced by autonomous vehicles (AVs), AVs' control strategy will be different from human driving behavior and their
traffic flow will be different from what we observe nowadays.
In this paper we would like to understand two questions:
\begin{enumerate}
  \item What is the new car-following control strategy of AVs at micro-scale?
  \item What is the new traffic flow theory for AVs at macro-scale?
\end{enumerate}
Human driving behavior has been extensively modeled at both micro- and macro-scales. At micro-scale,
 car following models (CFMs) treat each vehicle as a discrete entity with a constant length, whose dynamic location and velocity is computed from an underlying ordinary differential equation (ODE) system \cite{newell1961nonlinear,gipps1981behavioural,bando1995dynamical,brackstone1999car,zhang1999mathematical,zhang2005car,zheng2011applications,chen2012behavioral,orosz2010traffic,cui2017stabilizing,yang2010development,huang2017accelerated,zhao2017accelerated,zhao2018accelerated,liu2015safe,liu2016enabling,gong2016constrained,gong2018cooperative}.
CFMs assume local interactions among vehicles and local information from neighboring vehicles.
The modeling of each agent requires tracking and keeping records of surrounding agents.
Due to dynamic and volatile characteristics of traffic flow, the interacting agents and their topology may change quickly.
Real time design strategies may become extremely difficult to implement for heavy traffic scenarios,
as the associated microscopic control  mechanism may not be scalable for many vehicles. Moreover, it is also not easy to account for global traffic information obtained from vehicle connectivity.
In contrast, macroscopic traffic flow models treat one vehicle as a particle without occupying any space.
Traffic flow is then described by the continuum density distribution and velocity field
solved from partial differential equations (PDEs) \cite{newell1961nonlinear,newell1993simplified,lighthill1955kinematic,zhang1998theory,zhang2002non,jin2003distribution,daganzo2005variational,lebacque2005first}.

AVs, on the other hand, exhibit distinct driving behavior from HVs
and thus call for new models and theories for both microscopic vehicle control and macroscopic traffic flow.

\subsubsection{Microscopic Longitudinal Controller of AVs}
To prepare AVs to drive on public roads, safe and efficient controller design of autonomous driving is 
a top  priority.
AV controls can be categorized into
longitudinal control (i.e., the car-following scenario)
and lateral control (i.e., the lane-change scenario).
Longitudinal control has been studied in various scenarios, including:
platooning \cite{gong2016constrained,zhou2017rolling,wei2017dynamic,li2018nonlinear},
speed harmonization \cite{ma2016freeway,malikopoulos2018optimal,arefizadeh2018platooning},
longitudinal trajectory optimization \cite{wei2017dynamic,li2018piecewise},
and eco-approach and departure
at signalized intersections \cite{altan2017glidepath,hao2018eco,yao2018trajectory}.
Connected adaptive cruise control (CACC) is the most extensively studied longitudinal controller for AV platooning
\cite{swaroop1996string,darbha1999intelligent,
  ioannou1993autonomous,rajamani2001experimental,van2006impact,naus2010string,vanderwerf2001modeling,
  kesting2010enhanced,
  schakel2010effects,naus2010string,ploeg2011design,milanes2014cooperative,milanes2014modeling,talebpour2016influence,jin2014dynamics,wei2017dynamic,li2018piecewise}. 

This paper is primarily focused on AVs' \textbf{longitudinal velocity control} in the car-following scenario.
We formally formulate the problem as follows.
\begin{definition} (Problem Statement):
There are $N$ autonomous vehicles indexed by $i\in\{1,2,\dots,N\}$ driving in one direction on a closed highway without any entrance nor exit, with initial positions $x_{1,0},\dots,x_{N,0}$.
Each car aims to select its optimal velocity control by minimizing its driving cost functional pre-programmed by its manufacturer over the predefined planning horizon $\left[0, T\right]$.
We would like to investigate a scalable velocity control strategy for a large number of AVs.
\end{definition}
The modeling details will be discussed in next sections.
To develop microscopic AV controls, one needs to design autonomous driving behavior by some underlying dynamical models for AVs.
A majority of studies simply tailor AVs' behavior on that of HVs by tweaking behavioral parameters
(e.g., shorter reaction time or headway \cite{levin2016multiclass}
or k-vehicle ahead information \cite{jin2014dynamics}),
in which AVs are essentially human drivers but react faster, ``see" farther, and ``know" the road environment better.
The models proposed in those studies may not capture AVs' dynamic learning capabilities. Such learning capabilities are modeled by the model predictive control (MPC) or Stackelberg games in some studies \cite{wang2015game,Yu2018}. However, those studies suffer from scalability issues when the number of AVs becomes large.



\vspace{-6pt}
 \subsubsection{Macroscopic Traffic Flow Theory for AVs}
The connections between CFMs and macroscopic traffic flow models have been established either using continuum limit or via change of coordinates, states or formulations. In the former case, a macroscopic model is the limit of a CFM as the number of cars tends to infinity, which may be shown rigorously using theory from conservation laws \cite{RSMUP_2014__131__217_0,rossi2014justification} or measure theory \cite{di2015rigorous}. However, such continuum limit results are only known for very few CFMs while the general mathematical theories are still not available.  As an alternative, a macroscopic model can be transformed into different coordinates, different state variables or its variational formulation, so that consistency can be established between the transformed systems and specific CFMs.
For instance,
\cite{daganzo2006traffic,daganzo2006variational} showed that the Newell's CFM is a discrete form of the LWR model with the Greenshields fundamental diagram in Lagrangian coordinates.
Along this line, \cite{helbing2009derivation} studied a nonlocal second order model
and \cite{Laval2013} studied three representations of the LWR model using different coordinates and state variables.
Later, \cite{jin2016equivalence} established a more general framework bridging between a family of CFMs and macroscopic traffic flow models.

Similarly,
deriving a new traffic flow theory for pure AVs
requires the establishment of a micro-macro connection from AV control models.
There exist limited studies that characterized traffic flow theories from their respective microscopic controls of connected and automated vehicles using gas-kinetic theory \cite{ngoduy2009continuum,porfyri2015stability}.
In contrast, a majority of researchers simply derived new fundamental diagrams
using the existing traffic flow theory framework:
Assuming that AVs posses shorter reaction time in car-following, AVs' fundamental diagram has the same free-flow speed but steeper congestion curve \cite{levin2016multiclass} compared to HVs.


In this paper, we aim to derive a macroscopic game-theoretic model from AV microscopic longitudinal control
and fill two gaps in the existing literature by
(i) proposing an efficient computational framework for the longitudinal control of a platoon of AVs in the car-following scenario
and (ii) developing a traffic flow theory for AVs.

Traditional macroscopic traffic flow models are often classified into two categories: first-order models such as the Lighthill-Whitham-Richards (LWR) model \cite{lighthill1952sound,lighthill1955kinematic} and higher-order models such as the Payne-Whitham (PW) model \cite{Payne1971,Whitham2011} and the Aw-Rascle-Zhang (ARZ) model \cite{aw2000resurrection,zhang2002non}. The classification is based on different control variables. First-order models assume drivers control their speeds according to the traffic density while higher-order models prescribe a relationship between drivers' acceleration rates and the traffic density as well as drivers' speeds. The mean field game presented in this paper, which only models AVs' velocity controls, may be seen as the game-theoretic analogue to traditional first-order models.
 Similar to the extensions taken from traditional first-order models to higher-order models, one can incorporate more factors such as the acceleration rate and develop higher-order MFGs into our game-theoretic modeling framework. Table~\ref{tab:class} shows the classification of both traditional and game-theoretic macroscopic traffic flow models.
This paper primarily focuses on first-order MFGs and will leave the discussion on higher-order MFGs in future research.
\begin{table}[htbp]
  \centering
  \caption{Classification of macroscopic traffic flow models}
  \label{tab:class}
  \begin{tabular}{|c|c|c|}
    \hline
    \diagbox{Type}{Control} & Speed & Acceleration rate\\\hline
    Traditional & First-order (e.g., LWR) & Higher-order (e.g., PW/ARZ)\\\hline
    Game-theoretic & First-order MFGs & Higher-order MFGs\\\hline
  \end{tabular}
\end{table}


 \subsection{Literature review}\label{sec:lit}


Assuming connectivity between predecessors and followers as well as between platoon leaders and followers,
CACC contains two control policies:
constant spacing (CS) \cite{swaroop1996string,darbha1999intelligent,swaroop2001direct} and
constant time headway (CTH) \cite{ioannou1993autonomous,rajamani2001experimental,van2006impact,naus2010string,vanderwerf2001modeling,zhou2017rolling,arefizadeh2018platooning,stern2018dissipation}.
These two policies can be formulated as a linear time invariant system (LTI) \cite{swaroop1994comparision} with disturbances to dynamic and measurement dynamics \cite{zhou2017rolling} 
or a model predictive control (MPC) system
with distributed control \cite{wang2016cooperative,gong2016constrained,gong2018cooperative}.

AVs longitudinal acceleration control can also be modeled using nonlinear car following models (CFMs).
The most widely used CFMs for AVs are Intelligent Driver Model (IDM) \cite{treiber2000congested,kesting2008adaptive,kesting2010enhanced,schakel2010effects,naus2010string,ploeg2011design,milanes2014cooperative,milanes2014modeling,talebpour2016influence,cui2017stabilizing,wu2017emergent,talebpour2018effect,wu2018stabilizing,yang2010development,huang2017accelerated,zhao2017accelerated,zhao2018accelerated}
and Optimal Velocity Model (OVM) and its variants with heterogeneous communication delay or dynamic uncertainty 
\cite{qin2013digital,jin2014dynamics,qin2017scalable,jin2018connected,jin2018experimental}.
Unlike OVM, IDM takes safety into consideration and is thus collision-free.
All the aforementioned studies aim to develop a string stable car-following controller in order to smoothen traffic flow and prevent stop-and-go waves.
But none of them considers control and physical safety constraints \cite{gong2018cooperative}.
In other words, interactions among vehicles are not explicitly modeled \cite{li2018nonlinear}.

To address the above challenges, 
some researchers model a full penetration of AVs as a multi-agent system (MAS), wherein 
every AV interacts among one another through physical interactions in traffic.
A majority of studies that capture the interactions of vehicles assume that each vehicle 
carries out
 a sequence of accelerations over a finite time horizon
by optimizing a common or an individual objective function.
Vehicles interact 
among themselves
through the common or individual cost function as well as safety constraints.
Depending on the objective functional form,
these models can be further divided into two classes:
\emph{cooperative control} and \emph{non-cooperative game}.

Cooperative control has been widely studied in multi-robotic systems.
In light of multi-robotic-interaction, robots interact with one another and choose optimal policies by predicting others behavior.
Neighboring robots trajectories are treated as hard safety constraints or boundaries for robots motion planning.
Such modeling has been critical in multi-robot collision avoidance and human-robot interaction \cite{liu2015safe,liu2016enabling}.
A cooperative AV system is a multi-vehicle system that can be controlled to stabilize traffic flow and smoothen traffic jam \cite{wang2016cooperative,gong2016constrained,gong2018cooperative}, to optimize driving comfort \cite{wang2014rolling2,zhou2017rolling}, and to improve fuel efficiency \cite{wang2014rolling1,yao2018trajectory}. 
To reduce computational burdens, a distributed algorithm is usually designed and implemented on each vehicle \cite{wang2016cooperative,gong2016constrained,gong2018cooperative}.

Compared to cooperative AV control, the non-cooperative interactions among\\ \\ AVs are relatively understudied.
Game theory is a natural approach to model the non-cooperative strategic interactions among AVs,
in which each AV solved an MPC \cite{wang2015game} or MDP \cite{dreves2018generalized,li2018game,sadigh2016planning,lazar2018maximizing}.
In the game theoretic framework, cars are referred to as ``agents" or ``players".
\cite{wang2015game} formulated the discrete lane change and continuous acceleration selection of AVs as a differential game, where the agents' optimal strategies are obtained from solving optimal control problems.
\cite{Talebpour2015} modeled lane-changing behavior as a two-person non-zero-sum non-cooperative game under (in)complete information.
\cite{yoo2012stackelberg,yoo2013stackelberg,kim2014game,yu2018human} developed a Stackelberg game among multiple AVs in driving or merging and a mixed-motive game in lane-changing.
\cite{dreves2018generalized} modeled multiple AVs acceleration and steering angle velocity selection at intersections with the goal of avoiding collisions.  
The human-cyber-physical systems (h-CPS) community extends multi-agent systems to hybrid AVs interacting with human drivers. 
For example, \cite{sadigh2016planning} designed ``local interactions" between an AV and a human driver to drive efficiently and maximize road capacity, while \cite{lazar2018maximizing} generalized the model to several AVs and HVs.
\cite{li2018game} assumed that human drivers choose driving policies using hierarchical reasoning while AVs optimize car-following and lane-changing strategies based on a Stackelberg game.   
The outcomes of all the aforementioned game-theoretic models are equilibrium driving strategies.
The computation of equilibrium may become extremely challenging 
 when the number of coupled agents becomes large.
To get around, \cite{wang2015game} applied Model Predictive Control (MPC) instead of computing an equilibrium. 
\cite{dreves2018generalized} solved a generalized Nash equilibrium by summing up all vehicles objective functions, which is essentially a cooperative control.
\cite{lazar2018maximizing} assumed that AVs can directly 
perform
optimization based upon their predictions of human driver actions rather than human's actual strategies. 
None of these studies investigated quantitatively how close the approximate solutions are to the original differential games.  
Nevertheless, because the available game-based control algorithms suffer from scalability issues, all the aforementioned studies had to constrain their applications to a limited number of AVs. 
A scalable traffic simulation framework was developed where AVs learned their optimal driving policies in a multi-agent RL environment \cite{wu2017emergent,wu2017flow,wu2017framework},
but the trained policies suffered from a lack of interpretability.

%
%


All the aforementioned studies focus on AVs longitudinal or lateral controls in discrete games, which suffers from scalability issues.
Thus a scalable theory and algorithm applicable for a large number of coupled AV controllers is urgently needed.

\subsection{Contributions of this paper}

Instead of hypothesizing explicitly how AVs drive,
our goal is to design future AVs as rational, utility-optimizing agents
that play best driving strategies.
In other words,
AVs are intelligent agents programmed by manufacturers to minimize driving costs as a trade-off between traffic safety and efficiency.
Any deviation from driving with the best strategies will increase AVs' individual costs.

Game theory is a natural tool to model the equilibrium of interacting utility-optimizing agents.
Given that a large number of interacting AVs are designed to select velocity controls continuously,
we seek an innovative \emph{game-theoretic} tool,
i.e., the mean field game, 
for complex multi-agent dynamic modeling \cite{lasry2007mean,huang2006large}.
Mean field approximation allows for the translation of microscopic behaviors and interactions of agents to a macroscopic level. 
Most importantly,
we will show later in this paper that
MFG embodies classical traffic flow models with behavioral interpretation,
thereby providing a new traffic flow theory for AVs.\\

\section{The modeling framework} 

This paper contributes to the state-of-the-art of AV controller design by characterizing the interplay between discrete differential games and continuous mean field games.
Mean field game (MFG) is a game-theoretic model used to describe complex multi-agent dynamic systems \cite{lasry2007mean,huang2006large}. It has become increasingly popular in designing new decision-making processes for finance \cite{gueant2011mean,lachapelle2010computation},  engineering \cite{djehiche2016mean,couillet2012electrical}, social science \cite{degond2014large}, pedestrian crowds modeling \cite{lachapelle2011mean,burger2013mean} and traffic \cite{kachroo2016inverse,chevalier2015micro,huang2019stabilizing,huang2020scalable}.
MFG is a micro-macro model which allows one to define individuals on a microscopic level as rational utility-optimizing agents 
while translating rich microscopic behaviors to macroscopic models. 
The basic idea is to exploit the ``smoothing" effect of large numbers of interacting individuals. Instead of solving a long list of highly coupled equations that depict the interactions among different players, MFG assumes that each player only reacts to a ``mass" that results from an aggregate effect of all the players. Such an approach is called mean field approximation and helps to simplify the complex multi-agent dynamic systems on a macroscopic level.

This paper models AVs' rational and intelligent driving behavior under the mean field game framework. We shall make the following AVs' behavioral assumptions:
\begin{itemize}
  \item Each AV observes global in space traffic state information on the road.
  \item Each AV plans its velocity control in a time horizon by anticipating others' behaviors.
  \item AVs act to utilize their predefined driving costs over the time horizon in a non-cooperative way.
\end{itemize}

\begin{figure}[H]\centering
  \includegraphics[width=.7\textwidth, keepaspectratio=true]{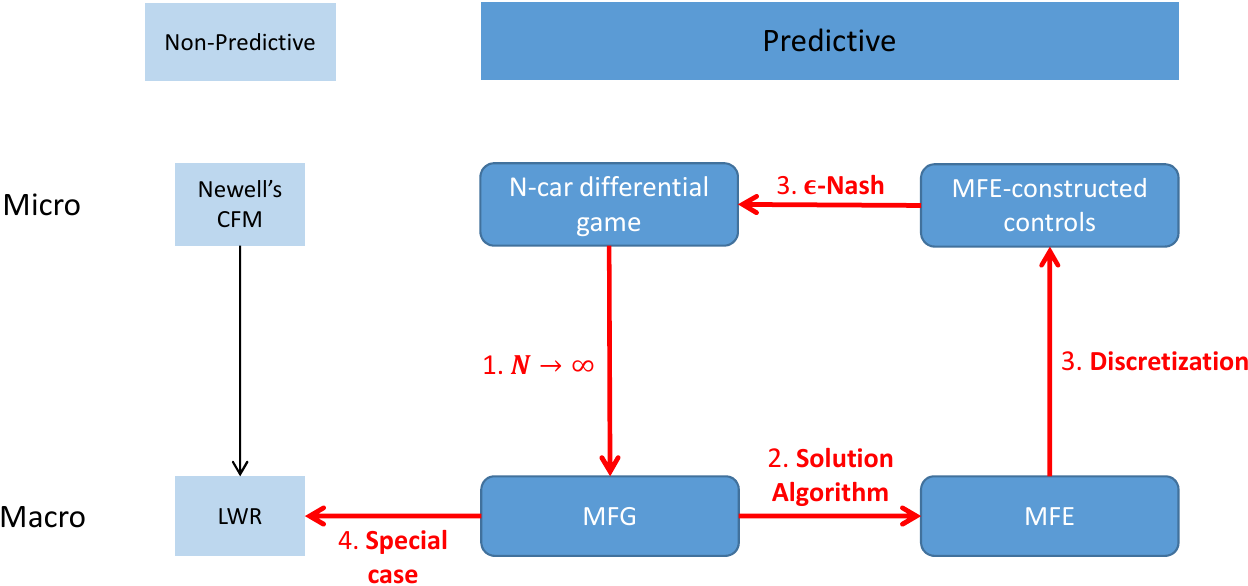}
  \vspace{-0.1in}
  \caption{From Micro to Macroscopic Traffic Flow Models}
  \label{fig:sche}
\end{figure}

Four major components of this paper are elaborated below  (as shown in Figure~\ref{fig:sche}):
\begin{enumerate}
  \item A mean field game is derived from the limiting differential game as the number of AVs tends to infinity. The mean field game is a coupled forward-backward PDE system that models AVs' non-cooperative velocity selections at a macroscopic scale.
  The existing research on the application of mean field games to transportation domain solely worked on specific objective functions \cite{chevalier2015micro,kachroo2016inverse}.
  In contrast, we systematically derive the forward continuity equation and the backward Hamilton-Jacobi-Bellman (HJB) equations with a family of more general objective functions using mean field approximation.

  \item An equilibrium solution, denoted by mean field equilibrium (MFE), is solved from the mean field game. AVs' optimal velocity control strategies are represented by the MFE at a macroscopic level.
  Existing algorithms for computing the MFE are mainly designed for a short planning horizon \cite{chevalier2015micro} or a special family of cost functions \cite{lachapelle2011mean,benamou2017variational}. In this paper we develop a new algorithm that works with a longer planning horizon and more general cost functions. Our algorithm is based on finite difference and multigrid preconditioned Newton's method.

  \item A tuple of AVs' discrete controls are constructed from the discretization of a continuous MFE.
  We test different numbers of AVs and different objective functions to illustrate the accuracy of MFE-constructed controls as an $\epsilon$-Nash equilibrium of the original differential game.
  The results show a consistent trend that the continuous equilibrium solution provides a good approximation to AVs' non-cooperative individual controls when the number of AVs is large.
  This construction method addresses the scalability issue faced by many existing literature \cite{wang2015game,dreves2018generalized,lazar2018maximizing}.


  \item The proposed mean field game can also be treated as a macroscopic traffic flow model. It models AVs' aggregated behavior assuming AVs are predictive and rational agents. Along this line, we first establish connections between the mean field game and the traditional LWR model rigorously. 
  Then we present some possible AV driving objective functions whose respective mean field games show interesting traffic patterns.
\end{enumerate}

%

The remainder of this paper is organized as follows. Section~\ref{sec:prob} introduces AVs' differential game as an extension to one AV's optimal longitudinal control problem. In Section~\ref{sec:MF}, the macroscopic MFG is derived from AVs' differential game with some assumptions. In Section~\ref{sec:flow}, we illustrate connections between MFG and the traditional LWR model in a general framework and present two MFG examples modeling AVs' kinetic energy, driving efficiency and safety. Then, Section~\ref{sec:sol} is devoted to a new algorithm to solve MFG numerically based on Newton's method. In Section~\ref{sec:micro-macro}, we construct a tuple of AVs' discrete controls from the continuous MFG solution and characterize their accuracy as an approximate equilibrium of the original differential game. Conclusions and future research directions follow in Section~\ref{sec:conclude}.

\section{From optimal control to differential game}\label{sec:prob}

We have seen a growing interest in applying optimal control theory to model AVs' predictive driving strategies in car-following and lane-change scenarios \cite{wang2014rolling1,wang2014rolling2,gong2016constrained,zhou2017rolling}.
In this section, we briefly introduce how to formulate a single AV's longitudinal control as an optimal control problem
and then extend it to a differential game among multiple AVs.

\subsection{Optimal longitudinal control of one car}\label{sec:oc_onecar}

Assume that there are $N$ AVs indexed by $i\in\{1,2,\dots,N\}$ driving in one direction on a closed highway of length $L$ without any entrance nor exit.
Denote the $i^{\text{th}}$ car's position at time $t$ by $x_i(t)$ and speed by $v_i(t)$. Fix a finite period of time $[0,T]$ where $T>0$, the cars' motions on $[0,T]$ are dictated by the following dynamical system:
\begin{align}
  \dot{x}_i(t)=v_i(t),\quad x_i(0)=x_{i,0},\quad i=1,2,\dots,N, \label{eq:cardynamics}
\end{align}
where,\\
$\dot{x}_i(t)$: the shorthand notation of $\frac{dx_i(t)}{dt}$;\\
$x_{i,0}$: the $i^{\text{th}}$ car's initial position at the beginning time $t=0$.\\
We use the notation $x_i(t)=x_i(t,v_i(\cdot),x_{i,0})$, $i=1,2,\dots,N$ for simplicity but keep in mind that $x_i(\cdot)$ depends on both $v_i(\cdot)$ and $x_{i,0}$.

For any $i=1,2,\dots,N$, suppose the $i^{\text{th}}$ car knows other cars' speeds:
\begin{align}
  \bm{v}_{-i}(t)=[v_1(t),\cdots,v_{i-1}(t),v_{i+1}(t),\cdots,v_N(t)]^T,
\end{align}
and positions:
\begin{align}
  \bm{x}_{-i}(t)=[x_1(t),\cdots,x_{i-1}(t),x_{i+1}(t),\cdots,x_N(t)]^T,
\end{align}
for $t\in[0,T]$. To select an optimal driving speed profile, the $i^{\text{th}}$ car solves an optimal control problem over the planning horizon $[0,T]$ .

Define the $i^{\text{th}}$ car's \emph{driving cost functional} as:
\begin{align}
  J_i^N(v_i,\bm{v}_{-i})=  \underbrace{\int_0^T \underbrace{f_i^N \left(v_i(t),x_i(t),\bm{x}_{-i}(t)\right)}_{\text{cost function}} \,dt}_{\text{running cost}}+\underbrace{V_T\left(x_i(T)\right)}_{\text{terminal cost}},\label{eq:optmctrlcost}
\end{align}
where,\\
$\int_0^T f_i^N \left(v_i(t),x_i(t),\bm{x}_{-i}(t)\right) \,dt$: the \emph{running cost} over the entire planning horizon;\\ 
$f_i^N(\cdot)$: the \emph{cost function} that quantifies driving objectives such as efficiency and safety;\\
$V_T(x_i(T))$: the \emph{terminal cost} representing the $i^{\text{th}}$ car's preference on its final position at time $T$.

We assume that all cars have the same free flow speed denoted by $u_{\text{max}}$. It is natural to require that the $i^{\text{th}}$ car's speed remains nonnegative 
and does not exceed $u_{\text{max}}$.  Mathematically, this means that
\begin{align}
  \mathcal{A}=\left\{v(\cdot):\ 0\leq v(t)\leq u_{\text{max}},\,\forall t\in[0,T]\right\},
\end{align}
is the admissible set of the $i^{\text{th}}$ car's speed selections. The $i^{\text{th}}$ car tries to obtain an optimal velocity control $v_i^*(\bm{v}_{-i}(\cdot),t)$ on the planning horizon $[0,T]$ such that:
\begin{align}
  J_i^N(v_i^*,\bm{v}_{-i})\leq J_i^N(v_i,\bm{v}_{-i}),\ \ \forall v_i\in\mathcal{A}.
\end{align}
$v_i^*(\bm{v}_{-i}(\cdot),t)$ depends on other cars' speeds $\bm{v}_{-i}(\cdot)$ through their trajectories $\bm{x}_{-i}(\cdot)$. We will use the notation $v_i^*(t)$ for simplicity.
When one car selects its own driving speed over the predefined planning horizon while everybody else does so simultaneously, a non-cooperative \emph{differential game} forms.

\subsection{\textit{N}-Car differential game}

Differential games can be regarded as extensions of non-cooperative Nash games in dynamic systems.
In a differential game, a finite number of players solve their individual optimal control problems while those optimal control problems are coupled through the dependency of one's cost functional on the others' actions \cite{basar1999dynamic}.
Along this line, we formulate the $N$-car differential game for AVs extending the one-car optimal control problem in Section~\ref{sec:oc_onecar}:

$N$ AVs indexed by $i\in\{1,2,\dots,N\}$ are driving in one direction on a closed highway without any entrance nor exit, with initial positions $x_{1,0},\dots,x_{N,0}$.
Each car aims to select its optimal velocity control by minimizing its driving cost functional defined in Eq.~\eqref{eq:optmctrlcost} over the predefined planning horizon $\left[0, T\right]$. A Nash equilibrium of the game is a tuple of controls $v^*_1(t),v^*_2(t),\dots,v^*_N(t)$ satisfying:
\begin{align}
  J_i^N(v^*_i,\bm{v}^*_{-i})\leq J_i^N(v_i,\bm{v}^*_{-i}),\ \ \forall v_i\in\mathcal{A},\quad i=1,\dots,N.\label{eq:DGE0}
\end{align}

It is generally difficult to solve an equilibrium when $N$ is large, because it involves solving $N$ coupled optimal control problems \cite{cardaliaguet2010notes}.
The goal of this paper is to develop a scalable framework to solve approximate equilibria for a family of $N$-car differential games by resorting to mean field approximation.

The underlying rationale of the developed methodology is articulated as follows:
\emph{
\begin{enumerate}
  \item Rather than solving the $N$-car differential game directly, we turn to its limit as the number of cars $N\to\infty$, i.e., a mean field game (Section~\ref{sec:MF}).
  \item A numerical algorithm is developed to solve the mean field game (Section~\ref{sec:sol}).
  \item The equilibrium solution of the mean field game is used to construct a tuple of discrete controls and those controls are verified to be an $\epsilon$-Nash equilibrium of the original $N$-car differential game
  by numerical examples (Section~\ref{sec:micro-macro}).
\end{enumerate}
}

\section{From differential game to mean field game}\label{sec:MF}


When the number of cars $N\to\infty$, one goes from the $N$-car differential game to a mean field game (MFG). The MFG is essentially a differential game with an infinite number of agents so that the interactions between any two individuals are ignorable. Instead, any individual reacts only to the ``mass'' of all agents. The ``mass'' then evolves with the aggregated behavior of all agents' motions. Two partial differential equations are developed to describe the MFG:
\begin{enumerate}
  \item A backward Hamilton-Jacobi-Bellman (HJB) equation: a generic car's speed selection is formulated as an optimal control problem where the generic car computes its driving cost associated to a cost function based on its prediction on the evolution of total ``mass''. The HJB equation is then derived from the optimal control problem. The solution of the HJB equation provides optimal costs and optimal velocity control strategies for all cars. The HJB equation is solved from $t=T$ to $t=0$ backward.
  \item A forward continuity equation: it is derived from the conservation of cars. The solution of the continuity equation describes the ``mass'' evolution arising from all cars' motions. The continuity equation is solved from $t=0$ to $t=T$ forward.
\end{enumerate}

\begin{figure}[H]\centering%
  \includegraphics[width=.85\textwidth, keepaspectratio=true]{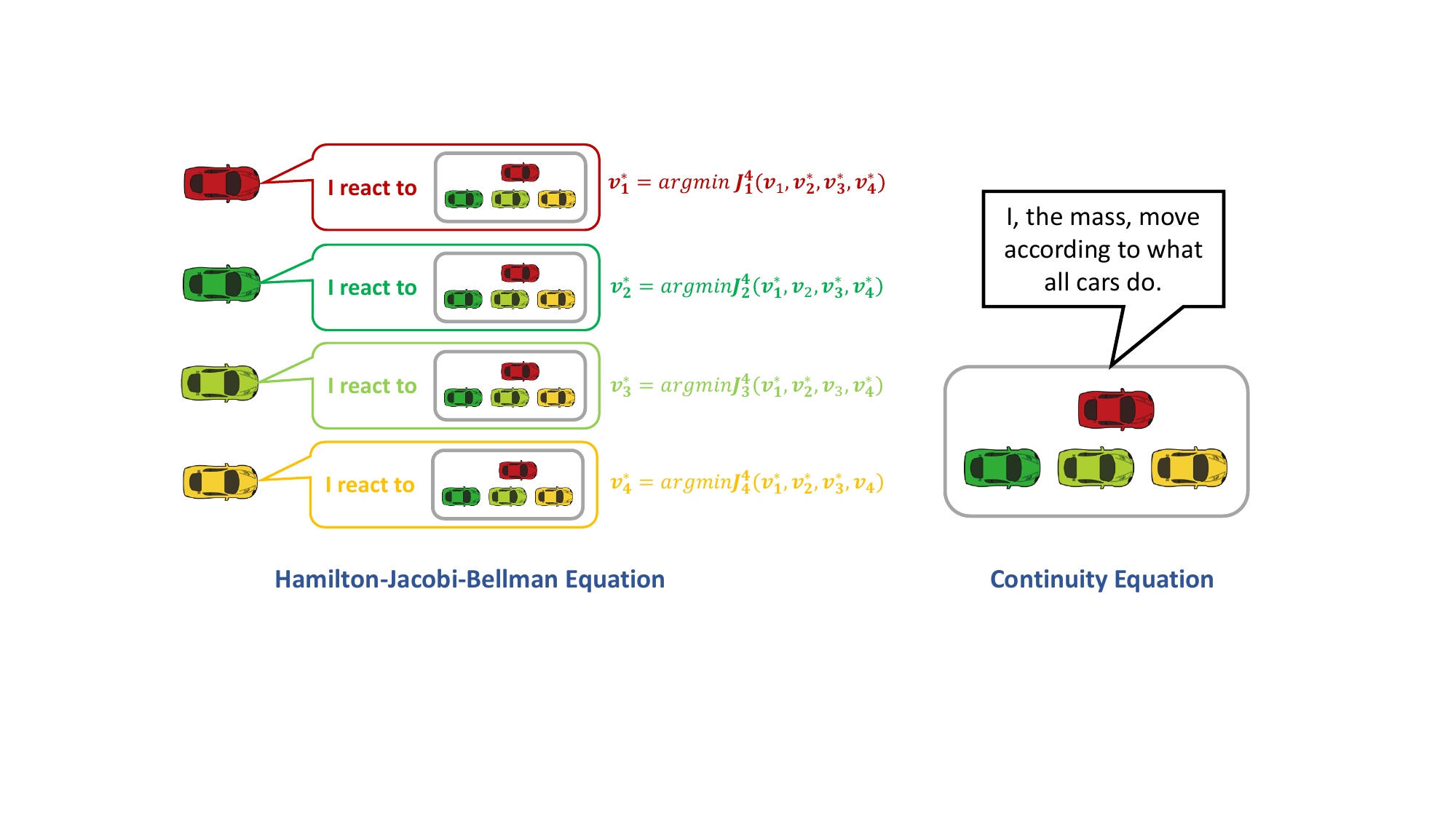}
  \caption{\small From an $N$-car differential game to MFG (adapted from \cite{MFGfish})}
  \label{fig:hjb_fpk}
\end{figure}

The MFG is a coupled system of the forward continuity equation and the backward HJB equation. At the mean field equilibrium, the total ``mass'' evolution coincides every car's prediction.
Figure~\ref{fig:hjb_fpk} shows a simple example of four cars to provide an intuitive explanation of these two equations.
 Each car is a rational agent aiming to minimize a driving cost defined in Eq.~\eqref{eq:optmctrlcost}, leading to a system of four coupled optimal control problems, one for each car.
 As $N$ goes large,
 the HJB equation can be derived from these coupled problems
 and the continuity equation can be derived from the trajectories of all cars.

In this section, we will formally derive the HJB and continuity equations from the $N$-car differential game using mean field approximation.

\subsection{Mean field limit}


The general idea of moving from the microscopic $N$-car differential game to the macroscopic MFG is to take a mean field limit by letting the number of cars in the system go to infinity.
To allow us to take the limit, we need to first make two homogeneity assumptions:
\begin{enumerate}
  \item[(A1)] All cars are indistinguishable.
  \item[(A2)] All cars have the same form of cost function.
\end{enumerate}

It should be mentioned that the above assumptions may be relaxed.
For example, (A2) can be relaxed if multi-class traffic is the subject of study. 
In this paper we mainly focus on single-class AVs and leave multi-class models to the future work.  

Provided that the $N$-car differential game satisfies the above assumptions, we will derive a MFG in four steps:
\begin{enumerate}
  \item We reformulate the driving cost functional defined in Eq.~\eqref{eq:optmctrlcost} by introducing a smooth density (Section~\ref{sec:smoothing});
  \item We derive a generic car's optimal control problem from the differential game by taking the mean field limit when $N\to\infty$ (Section~\ref{sec:generic});
  \item We derive a set of HJB equations from the generic car's optimal control problem (Section~\ref{sec:hjb});
  \item We obtain an evolution equation and show that it is exactly the continuity equation widely used in macroscopic traffic flow models (Section~\ref{sec:evol}).
\end{enumerate}

\subsubsection{Step 1: Traffic Density Smoothing}\label{sec:smoothing}

Traffic density is a crucial quantity to manifest the macroscopic aspect of traffic flow.
Assumption (A1) enables us to replace states of individual cars in the driving cost functional defined in Eq.~\eqref{eq:optmctrlcost} by an aggregated traffic density.

More precisely, for any $i=1,\dots,N$, assumption (A1) implies that the function $f_i^N \left(v_i(t),x_i(t),\bm{x}_{-i}(t)\right)$ does not depend on the permutation of $x_1(t),x_2(t),\dots,x_N\break (t)$. According to \cite{cardaliaguet2010notes}, we can replace the arguments $x_i(t),\bm{x}_{-i}(t)$ in $f_i^N \left(v_i(t),x_i(t),\right.$\break $\left.\bm{x}_{-i}(t)\right)$ by an \emph{empirical density distribution} of $x_1(t),x_2(t),\dots,x_N(t)$, which is defined as:
\begin{align}
  \rho^N(x,t)=\frac1N\sum_{j=1}^N \delta\left(x-x_j(t)\right),
\end{align}
where $\delta(\cdot)$ is the Dirac mass.

However, $\rho^N$ is not a smooth function, leading to non-smoothness of the new driving cost functional. 
To resolve this issue, we first approximate $\rho^N$ using a smoothing kernel.

Suppose that $\xi(x)$ is a smoothing kernel which is smooth and nonnegative, and satisfies $\int_{\mathbb{R}}\xi(x)\,dx=1$. We take a smoothing parameter $\sigma>0$ and define the scaled kernel $\xi_\sigma(x)=\frac1\sigma \xi(\frac{x}{\sigma})$. The physical meaning of using the scaled kernel $\xi_\sigma(x)$ is that the $i^{\text{th}}$ car contributes to the density in a ``window" $[x_i(t)-\sigma,x_i(t)+\sigma]$ rather than only at the point $x_i(t)$ $(i=1,\dots,N)$ so that the density changes smoothly with location $x$. The smooth density distribution is defined as:
\begin{align}
  \rho^N_{\sigma}(x,t)=\frac{1}{N}\sum_{j=1}^N\xi_{\sigma}\left(x-x_j(t)\right).\label{eq:smooth}
\end{align}

With the smooth density, the $i^{\text{th}}$ car's cost function is rewritten as: 
\begin{align}
  f_i^N\left(v_i(t),x_i(t),\bm{x}_{-i}(t)\right) \triangleq f_i\left(v_i(t),\rho^N_\sigma(\cdot,t)\right),
\end{align}
where $v_i(t)$ is the car's speed and $\rho^N_\sigma(\cdot,t)$ is the traffic density over the whole road at time $t$. Generally $f_i$ may have arbitrary dependence on $\rho^N_\sigma(\cdot,t)$ as well as its spatial derivatives. To simplify, we make the assumption:
\begin{enumerate}
  \item[(A3)] The cost function only depends on the traffic density at the car's position.
\end{enumerate}

By assumption (A3) we can write:
\begin{align}
  f_i\left(v_i(t),\rho^N_\sigma(\cdot,t)\right) \triangleq f_i\left(v_i(t),\rho^N_\sigma(x_i(t),t)\right),
\end{align}
where $f_i(\cdot,\cdot)$ is a bivariate function of speed and density, $i=1,2,\dots,N$.

By assumption (A2), we have $f_1=f_2=\cdots=f_N=f$, where $f(\cdot,\cdot)$ is a bivariate cost function shared by all cars.
In summary, the $i^{\text{th}}$ car's driving cost becomes:
\begin{align}{}
  J_i^N(v_i,\bm{v}_{-i})=\int_0^T f\left(v_i(t),\rho^N_\sigma(x_i(t),t)\right)\,dt+V_T\left(x_i(T)\right). \label{eq:optmctrlcostmf}
\end{align}

It should be noted that the density information ahead of and behind the $i^{\text{th}}$ car is asymmetric. At any time $t_1\in[0,T)$, the car anticipates the model predicted density $\rho_\sigma^N(x_i(t_2),t_2)$ for a later time $t_2\in (t_1,T]$ to select its driving speed at time $t_1$. Since the $i^{\text{th}}$ car drives at positive speeds, it always holds that $x_i(t_2)\geq x_i(t_1)$. The fact yields that the cars ahead of the $i^{\text{th}}$ car may contribute to the density $\rho_\sigma^N(x_i(t_2),t_2)$ but never will the cars behind the $i^{\text{th}}$ car do so. Consequently, the $i^{\text{th}}$ car's velocity control is not influenced by the cars behind it.

\begin{definition}\label{def:dg}
  \begin{enumerate}
  \item \textbf{$N$-car mean field type differential game}  [DG]:

  $N$ AVs indexed by $i\in\{1,2,\dots,N\}$ are driving in one direction on a closed highway of length $L$ without any entrance nor exit, with initial positions $x_{1,0},\dots,x_{N,0}$.
  Each car aims to select its optimal velocity control by minimizing its driving cost functional defined in Eq.~\eqref{eq:optmctrlcostmf} over the predefined planning horizon $\left[0, T\right]$.

  \item \textbf{$N$-car mean field type differential game equilibrium}  [DGE]:

  A Nash equilibrium of the $N$-car mean field type differential game is a tuple of controls $v^*_1(t),v^*_2(t),\dots,v^*_N(t)$ satisfying:
  \begin{align}
    J_i^N(v^*_i,\bm{v}^*_{-i})\leq J_i^N(v_i,\bm{v}^*_{-i}),\ \ \forall v_i\in\mathcal{A},\quad i=1,\dots,N.\label{eq:DGE}
  \end{align}
  At equilibrium, no car can improve its driving cost by unilaterally switching its velocity control.
  \end{enumerate}
\end{definition}

We see from Eq.~\eqref{eq:optmctrlcostmf} that each car only responds to and contributes to the density distribution $\rho^N_\sigma$ of all cars through driving costs. Such a property allows us to take the \emph{mean field limit} of the game as $N$ tends to infinity.

\subsubsection{Step 2: Optimal Control of a Generic Car}\label{sec:generic}
We take the mean field limit in the following way: fix the ratio $L/N$, let $N\to\infty$ and $\sigma/L\to0$. Intuitively that means we fix the space headway and shrink the ``window" so that in the limiting case one car only sees a local density. Under the limit, using mean field approximation we replace $\rho_\sigma^N(x,t)$ that is computed from $N$ cars' positions by a continuum density distribution $\rho(x,t)$. Note that all cars are anonymous, we can ignore the index $i$ and consider a generic car starting from $x_0$ at $t=0$. Denote the car's velocity control by $v(t)$ and trajectory by $x(t)$ for $t\in[0,T]$, we rewrite Eq.~\eqref{eq:optmctrlcostmf} as:
\begin{align}
  J(v)=\int_0^T f\left(v(t),\rho(x(t),t)\right)\,dt+V_T\left(x(T)\right),\label{eq:genericarcost}
\end{align}
where its dynamic motion is described by
\begin{align}
  \dot{x}(t)=v(t),\quad x(0)=x_0,\label{eq:genericardyn}
\end{align}
and its velocity control $v(\cdot)$ is constrained by
\begin{align}
  0\leq v(t)\leq u_{\text{max}},\ \ \forall t\in[0,T].\label{eq:genericons}
\end{align}

\subsubsection{Step 3: HJB Equation}\label{sec:hjb}

{
 The generic car solves the optimal control problem Eqs.~\eqref{eq:genericarcost}\eqref{eq:genericardyn}\eqref{eq:genericons} to obtain its optimal velocity control $v^*(t)$, which depends on the generic car's initial position $x_0$. The initial position $x_0$ can be any position on the road.
Rather than solving an infinite number of optimal control problems for every initial position, we use dynamic programming and derive a set of HJB equations that characterize the optimality condition of the velocity control. Such an approach is widely used in optimal control theory \cite{bardi2008optimal}.

We first introduce the Bellman value function $V(x,t)$ and the optimal velocity field $u(x,t)$. $V(x,t)$ is defined as the optimal cost for the generic car starting from location $x$ at time $t$:
\begin{subequations}
\begin{align}
  V(x,t)&=\min\nolimits_{v:[t,T]\to [0,u_{\text{max}}]}\left[\int_t^T f\left(v(s),\rho(x(s),s)\right)\,ds+V_T\left(x(T)\right)\right],\label{eq:dp1}\\
  \text{s.t. }\dot{x}(s)&=v(s),\quad x(t)=x,\label{eq:dp2}
\end{align}
\end{subequations}
and $u(x,t)$ is defined as the car's speed at location $x$ and time $t$ when choosing the optimal control of Eqs.~\eqref{eq:dp1}\eqref{eq:dp2}.

From another point of view, $v^*(t)$ is the Lagrangian optimal velocity control while $u(x,t)$ is the Eulerian optimal velocity field. Once $u(x,t)$ is solved for all $x$ and $t$, the optimal cost of the original problem Eqs.~\eqref{eq:genericarcost}\eqref{eq:genericardyn}\eqref{eq:genericons} is given by $V(x_0,0)$ and the optimal control $v^*(t)$ is given by the feedback law:
\begin{align}
  v^*(t)&=u(x^*(t),t),\label{eq:feed1}\\
  \dot{x}^*(t)&=v^*(t),\quad x^*(0)=x_0.\label{eq:feed2}
\end{align}
}

Then we derive the HJB equations for $V(x,t)$ and $u(x,t)$ from Eqs.~\eqref{eq:dp1}\eqref{eq:dp2}. Suppose the generic car starts from position $x$ at time $t$.
Consider a small time step $\Delta t$, we can divide the driving cost in Eq.~\eqref{eq:dp1} into two parts on $[t,t+\Delta t]$ and $[t+\Delta t,T]$:
\begin{align}
  &\int_t^T f\left(v(s),\rho(x(s),s)\right)\,ds+V_T\left(x(T)\right)\notag\\
  =&\int_{t}^{t+\Delta t} f(v(s),\rho(x(s),s))\,ds+\int_{t+\Delta t}^T f(v(s),\rho(x(s),s))\,ds+V_T\left(x(T)\right).
\end{align}
Correspondingly, the generic car's decision process is also divided into two stages. First it selects the speed $v(t)=\alpha\in[0,u_{\text{max}}]$ on the horizon $[t,t+\Delta t]$. Then it moves to $x+\alpha\Delta t$ at time $t+\Delta t$ and selects its speed profile over the rest of the planning horizon $[t+\Delta t,T]$.

The running cost on $[t,t+\Delta t]$ is approximated by
\begin{align}
  \int_{t}^{t+\Delta t} f(v(s),\rho(x(s),s))\,ds=f(\alpha,\rho(x,t))\Delta t+O(\Delta t^2).
\end{align}
Note that from the new position $x+\alpha\Delta t$, the optimal cost on $[t+\Delta t,T]$ the car can obtain is $V(x+\alpha\Delta t,t+\Delta t)$.
By dynamic programming principle we have:
\begin{align}
  V(x,t)=\min\nolimits_{0\leq\alpha\leq u_{\text{max}}}\left\{f(\alpha,\rho(x,t))\Delta t+V(x+\alpha\Delta t, t+\Delta t)+O(\Delta t^2)\right\}.\label{eq:hjbderive1}
\end{align}
Take the first order Taylor's expansion of $V(x+\alpha\Delta t, t+\Delta t)$ near $(x,t)$, denote $V_t$ and $V_x$ the partial derivatives $\frac{\partial V}{\partial t}$ and $\frac{\partial V}{\partial x}$. Eq.~\eqref{eq:hjbderive1} yields:
\begin{align}
  V(x,t)=\min\nolimits_{0\leq\alpha\leq u_{\text{max}}}&\left\{f(\alpha,\rho(x,t))\Delta t+V(x,t)+\alpha\Delta tV_x(x,t)+\Delta tV_t(x,t)\right.\notag\\
  &+\left.O(\Delta t^2)\}\right..
\end{align}
Eliminating $V(x,t)$ from both sides, dividing both sides by $\Delta t$ and letting $\Delta t\to0$, we get:
\begin{align}
  V_t+\min\nolimits_{0\leq\alpha\leq u_{\text{max}}}\{f(\alpha,\rho)+\alpha V_x\}=0.\label{eq:appenhjb1s}
\end{align}
We assume that $f$ is strictly convex with respect to its first argument, the driving speed. Then we can introduce:
\begin{align}
  f^*(p,\rho)=\min\nolimits_{0\leq\alpha\leq u_{\text{max}}}\{f(\alpha,\rho)+\alpha p\},\ \forall p\in\mathbb{R},\label{eq:legendre}
\end{align}
so that $-f^*(-\cdot,\rho)$ is the Legendre transformation of $f(\cdot,\rho)$ for any $\rho$. Using $f^*$, Eq.~\eqref{eq:appenhjb1s} can be rewritten as:
\begin{align}
  V_t+f^*(V_x,\rho)=0.\label{eq:appenhjb1}
\end{align}
{
The strict convexity of $f$ with respect to speed yields the uniqueness of the minimizer in Eq.~\eqref{eq:legendre} that is given by $f^*_p(p,\rho)$ for any $p\in\mathbb{R}$, where $f_p$ denotes $f$'s derivative with respect to $p$.}
As a result, the optimal velocity field $u(x,t)$ is given by:
\begin{align}
  u=\argmin\nolimits_{0\leq\alpha\leq u_{\text{max}}}\{f(\alpha,\rho)+\alpha V_x\}=f^*_p(V_x,\rho).\label{eq:appenhjb2}
\end{align}
{
We highlight the convexity assumption on $f$ because
the car's optimal speed selection is not unique even in a single time step without the assumption. In real applications, it is reasonable to assume AVs’ utility satisfies the law of diminishing marginal returns \cite{shephard1974law,imprialou2016re}, i.e., increasing speed results in smaller increase in utility. As a corollary to the law, the utility should be concave with respect to speed. Then the convexity assumption follows from the fact that AVs' driving cost is just the negative of the utility.}

When $t=T$, Eq.~\eqref{eq:dp1} becomes $V(x,T)=V_T(x)$, which gives the terminal condition of the HJB equations.

Summarizing all above, given the density distribution $\rho(x,t)$, Eqs.~\eqref{eq:dp1}\eqref{eq:dp2} lead to the following HJB equations:
\begin{subequations}
\begin{numcases}{\mbox{(HJB)}\quad}
V_t+f^*(V_x,\rho)=0,\label{eq:hjb1}\\
u=f^*_p(V_x,\rho).\label{eq:hjb2}
\end{numcases}
\end{subequations}
with $V(x,T)=V_T(x)$. $V(x,t)$ and $u(x,t)$ are solved backward from the HJB equations.

\subsubsection{Step 4: Continuity Equation}\label{sec:evol}

When all cars follow the optimal velocity control, the aggregated density distribution $\rho(x,t)$ evolves according to the optimal velocity field $u(x,t)$ obtained from the HJB equations. An evolution equation can be derived from the conservation of cars:
\begin{align}
  \mbox{(CE)}\quad\rho_t+(\rho u)_x=0, \label{eq:continuity}
\end{align}
to describe the evolution of density $\rho(x,t)$ from some initial density distribution $\rho(x,0)=\rho_0(x)$. Eq.~\eqref{eq:continuity} is exactly the continuity equation (CE) widely used in traffic flow models \cite{orosz2010traffic}. Given velocity field $u(x,t)$ known, Eq.~\eqref{eq:continuity} is solved forward.

\subsection{Mean field game system}

Summarizing Eqs.~\eqref{eq:continuity}\eqref{eq:hjb1}\eqref{eq:hjb2}, when the HJB and continuity equations are coupled, we have the following MFG system with the cost function $f(u,\rho)$:
\begin{subnumcases}{\text{[MFG]}\quad}
  \mbox{(CE)}& $\rho_t+(\rho u)_x=0$,\label{mfg_cont}\\
  \mbox{(HJB)} & $V_t+f^*(V_x,\rho)=0$,\label{eq:mfg_hjb1}\\
   & $u=f^*_p(V_x,\rho)$.\label{eq:mfg_hjb2}
\end{subnumcases}



The associated initial and terminal conditions are provided by the initial density $\rho(x,0)=\rho_0(x)$ and the terminal cost $V(x,T)=V_T(x)$, respectively. The choice of boundary conditions depends on the traffic scenario. When cars drive on a ring road without any entrance nor exit, 
periodic boundary conditions are specified as: $\rho(0,t)=\rho(L,t),V(0,t)=V(L,t)$; When the road has an entrance at $x=0$ and an exit at $x=L$, we should impose the boundary conditions $\rho(0,t)=\rho_{\text{entr}}(t)$ representing the inflow at the entrance and $V(L,t)=V_{\text{exit}}(t)$ representing the boundary cost when cars leave the road at the exit. This paper will focus on the periodic boundary conditions.

Denote the system's solution by $\rho^*(x,t)$ and $u^*(x,t)$.
The optimal velocity field $u^*(x,t)$ is our primary focus and will thus be referred as the mean field equilibrium (MFE) in the subsequent analysis.

\begin{remark}
    \,[MFG] is the general MFG system with any cost function $f(u,\rho)$ that is strictly convex with respect to $u$. The existence and uniqueness of MFE for the general system remains to be investigated. MFGs with some special cost functions are shown to have a unique MFE, see discussions in Section~\ref{sec:mfgsep}.
 \end{remark}
\begin{remark}   
The MFG system derived here is usually called a
    \emph{non-viscous} MFG system, because we assume no stochasticity on cars' dynamics.
    Accordingly, [MFG] has no viscous terms such as $\rho_{xx}$ and $V_{xx}$. For theory on non-viscous MFG, 
    we refer to \cite{cardaliaguet2010notes,cardaliaguet2015weak}.
\end{remark}


\section{Mean field games in traffic flow}\label{sec:flow}

MFG shares the same continuity equation with traditional traffic flow models but characterizes cars' reactions to traffic congestions in a different way.
Traditional traffic flow models prescribe a relationship between traffic density and the car's speed or acceleration, while MFG models the car's speed selection as an optimal control problem with a prescribed cost function.
Based upon such understanding, MFG can be seen as a macroscopic traffic flow model that models AVs' predictive and rational driving behavior.
In this section we first establish connections between MFG and the traditional LWR model, and then present MFG examples by choosing appropriate cost functions to quantify AVs' driving objectives.

\subsection{Connections between MFG and LWR}
The Lighthill-Whitham-Richards (LWR) model \cite{lighthill1955kinematic,richards1956shock} is a representative of traditional traffic flow models, so it would be helpful to establish connections between MFG and LWR. \cite{kachroo2016inverse,chevalier2015micro} have revealed such connections focusing on a specific class of LWR with the Greenshields fundamental diagram. \cite{kachroo2016inverse} presented a cost function whose corresponding MFG takes the Greenshields LWR as a solution. \cite{chevalier2015micro} claimed that the Greenshields LWR is essentially a MFG with a specific cost function when drivers minimize their driving costs by selecting their driving speeds myopically .

\begin{figure}[H]
  \centering
  \includegraphics[width=.9\textwidth]{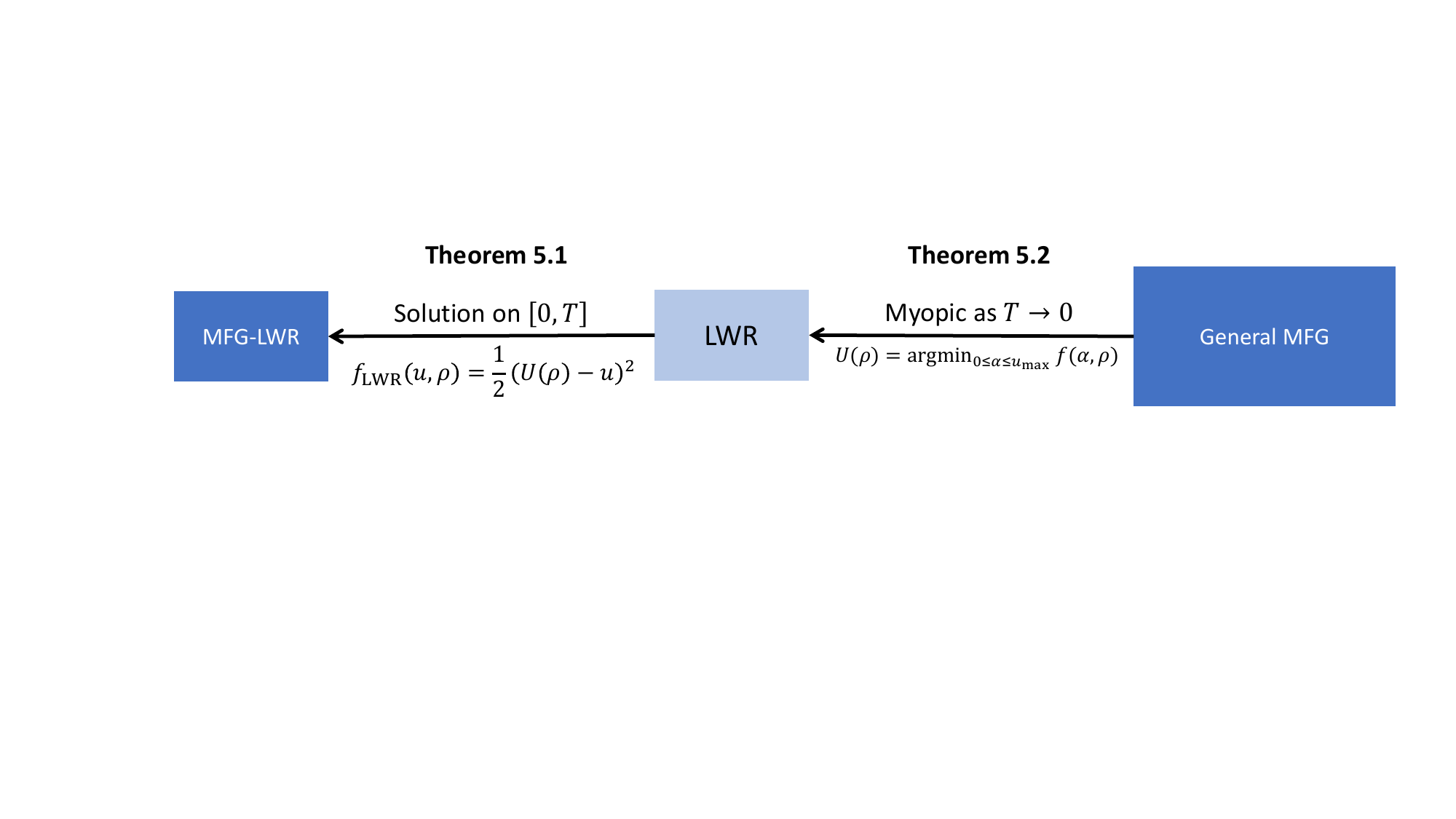}
  \caption{Connections between MFG and LWR}
  \label{fig:LWR}
\end{figure}

In the subsequent analysis we will establish connections between MFG and LWR from two perspectives, as shown in Figure~\ref{fig:LWR}.
(i) We coin a cost function for an arbitrary fundamental diagram and show that the LWR is a solution of the corresponding MFG;
(ii) LWR can also be seen as the myopic limit of MFG by letting the length of the planning horizon tend to zero, with a general family of cost functions.


\subsubsection{LWR as a solution to MFG}\label{sec:lwrsolmfg} 

Let us choose an arbitrary desired speed function $U(\rho)$. The corresponding LWR model is:
\begin{subequations}
\begin{numcases}{\text{[LWR]}\quad}
    \rho_t+(\rho u)_x=0,\label{eq:lwr_cont}\\
    u=U(\rho).\label{eq:lwr_speed}
\end{numcases}
\end{subequations}

Now we directly set the driving objective to be maintaining the LWR speed. There are infinite choices of respective cost functions. Here we artificially choose the following cost function:
\begin{align}
  f_{\text{LWR}}(u,\rho)=\frac{1}{2}(U(\rho)-u)^2.\label{eq:costfun_lwr}
\end{align}

Eq.~\eqref{eq:costfun_lwr} can be interpreted as quantifying the difference between the car's actual speed and desired speed. In other words, the car's objective is to keep not too far from human driving. The other reason to choose the cost function $f_{\text{LWR}}$ is that it relates to another cost function $f_{\text{NonSep}}$ modeling driving efficiency and safety, which will be shown later.

The cost function $f_{\text{LWR}}$ corresponds to the following MFG system:
\begin{subequations}
\begin{numcases}{\text{[MFG-LWR]}\quad}
  \rho_t+(\rho u)_x=0,\label{eq:mfg_lwr_cont}\\
  V_t+U(\rho)V_x-\frac12V_x^2=0,\label{eq:mfg_lwr_hjb1}\\
  u=U(\rho)-V_x.\label{eq:mfg_lwr_hjb2}
\end{numcases}
\end{subequations}

\begin{theorem}
  \label{thm:mfglwr}
  The solution of [LWR] is a solution of [MFG-LWR] under the conditions that: (i) [MFG-LWR] and [LWR] have the same initial density $\rho_0(x)$ and periodic boundary conditions; (ii) $V_T(x)=C$ where $C$ is an arbitrary constant for [MFG-LWR].
\end{theorem}
\begin{proof}
  Denote $\rho^*(x,t)$ and $u^*(x,t)$ the solution of [LWR]. Note that Eq.~\eqref{eq:lwr_cont} is the same as Eq.~\eqref{eq:mfg_lwr_cont}, it suffices to show that $\rho^*$ and $u^*$ satisfy the HJB equations~\eqref{eq:mfg_lwr_hjb1}\eqref{eq:mfg_lwr_hjb2} for some $V^*$. Take $V^*\equiv C$, then the terminal condition $V^*(x,T)=V_T(x)=C$ is satisfied and Eqs.~\eqref{eq:mfg_lwr_hjb1}\eqref{eq:mfg_lwr_hjb2} become a single equation $u^*=U(\rho^*)$, which is true from Eq.~\eqref{eq:lwr_speed}. So $\rho^*$, $u^*$ and $V^*\equiv C$ is a solution of [MFG-LWR].
\end{proof}

\begin{remark}
  The solution uniqueness of [MFG-LWR] is not studied in this paper. Proving the solution existence and uniqueness of MFG systems with general \emph{nonseparable} cost functions is mathematically challenging, such existence and uniqueness results are only obtained with short time horizon or small initial density \cite{ambrose2018existence}. Here we show the solution existence of [MFG-LWR] from Theorem~\ref{thm:mfglwr}. When the solution uniqueness holds for [MFG-LWR], its unique solution is the solution of [LWR] and the two systems [MFG-LWR] and [LWR] are equivalent. The equivalence between [MFG-LWR] and [LWR] is also supported by the numerical experiment. The rigorous proof will be left for the future research.
\end{remark}

\begin{remark}
  $V_T(x)=C$ means that cars have no preference on their final positions. One can specify the preference by imposing a non-constant terminal cost \cite{lachapelle2011mean}. In this paper we will always assume the terminal cost $V_T(x)=C$.
\end{remark}

Theorem~\ref{thm:mfglwr} will be verified with the Greenshields desired speed function:
\begin{align}
  U(\rho)=u_{\text{max}}(1-\frac{\rho}{\rho_{\text{jam}}}),\label{eq:greenshields}
\end{align}
later in the numerical experiment, where $\rho_{\text{jam}}$ is the jam density.

\subsubsection{LWR as the myopic limit of MFG}

To demonstrate the other connection between MFG and LWR, we consider a general cost function $f(u,\rho)$ and its corresponding MFG system [MFG].

Given the planning horizon $[0,T]$, a generic car selects its optimal velocity control to minimize the driving cost functional defined in Eq.~\eqref{eq:genericarcost}. If the generic car is myopic and does not concern the future, intuitively it will select the speed $u$ to minimize the instantaneous cost, i.e., $u=\argmin\nolimits_{0\leq\alpha\leq u_{\text{max}}} f(\alpha,\rho)$ at any time $t$, which leads to a LWR model with the desired speed:
\begin{align}
U(\rho)=\argmin\nolimits_{0\leq\alpha\leq u_{\text{max}}} f(\alpha,\rho)=f^*_p(0,\rho),\label{eq:fd}
\end{align}
according to Eq.~\eqref{eq:appenhjb2}.

{
To give a rigorous description of the myopic behavior, we define the myopic limit to be the limiting process when the length of the planning horizon $T\to0$. In the myopic limit, the anticipation effect of future traffic tends to zero. It is expected that the solution of the MFG will converge to the solution of the LWR with the desired speed function $U(\rho)$ defined in Eq.~\eqref{eq:fd}.}

\begin{theorem}\label{thm:myopic}
  Under the conditions that: (i) $f(u,\rho)$ is continuously differentiable, strictly convex with respect to $u$; (ii) the terminal cost $V_T(x)=C$ where $C$ is an arbitrary constant for [MFG]; (iii) there exists $T_0>0$ such that whenever $0<T\leq T_0$, with initial density $\rho_0(x)$ and periodic boundary conditions, [MFG] has a unique solution $\rho^{(T)}(x,t)$, $u^{(T)}(x,t)$ and $V^{(T)}(x,t)$ which are uniformly bounded up to second order derivatives on $0\leq x\leq L$, $0\leq t\leq T\leq T_0$. When $T\to0$ we have:
{
\begin{align}
  \lim_{T\to0}u^{(T)}(x,0)=U(\rho_0(x)),\ \forall x\in[0,L].
\end{align}
}
\end{theorem}
\begin{proof}
  There exists a constant $M>0$ such that $|V^{(T)}_{xt}(x,t)|\leq M$ for all $0\leq x\leq L$ and $0\leq t\leq T\leq T_0$. Integrate the inequality from $t=T$ to $t=0$, note that $V_x^{(T)}(x,T)=\frac{dV_T(x)}{dx}=0$ for all $0\leq x\leq L$, we get:
  \begin{align}
    |V^{(T)}_x(x,0)|\leq MT,
  \end{align}
  for all $0\leq x\leq L$ and $0\leq T\leq T_0$. Hence $V^{(T)}_x(x,0)\to0$ when $T\to0$, $\forall x \in[0,L]$.

  Since $f(u,\rho)$ is continuously differentiable and strictly convex with respect to $u$, $f_p^*(p,\rho)$ is continuous with respect to $p$. From Eq.~\eqref{eq:mfg_hjb2} we deduce that when $T\to0$:
  \begin{align}
    u^{(T)}(x,0)=f_p^*(V_x^{(T)}(x,0),\rho_0(x))\to f_p^*(0,\rho_0(x))=U(\rho_0(x)),\ \forall x\in[0,L].
  \end{align}
\vspace{-10pt}\end{proof}
\begin{remark}
  Typically a desired speed function $U(\rho)$ is supposed to satisfy certain conditions. For example: (i) $U'(\rho)\leq0$; (ii) $U(0)=u_{\text{max}}$; (iii) $U(\rho_{\text{jam}})=0$. In Theorem~\ref{thm:myopic}, $U(\rho)$ is computed from Eq.~\eqref{eq:fd}. The conditions on $U(\rho)$ are rewritten as: (i) $f^*_{p\rho}(0,\rho)\leq0$; (ii)$f^*_p(0,0)=u_{\text{max}}$; (iii) $f^*_p(0,\rho_{\text{jam}})=0$. Here the subscripts represent respective partial derivatives. Using the identity $f_u(f^*_p(p,\rho),\rho)=p$ between $f$ and $f^*$ \cite{rockafellar2015convex}, we can translate the conditions on $f^*$ to those on $f$. As a result, we require the cost function $f(u,\rho)$ to satisfy: (i) $f_{u\rho}(U(\rho),\rho)\geq0$; (ii) $f_u(u_{\text{max}},0)=0$; (iii) $f_u(0,\rho_{\text{jam}})=0$.  These conditions provide a way to calibrate the cost function from its myopic behavior.
\end{remark}
We can now interpret LWR from the perspective of MFG, which provides a richer behavioral foundation
and a more general and flexible framework.

\subsection{MFG examples}

As a micro-macro game-theoretic model, MFG can capture richer driving behaviors than LWR by choosing various cost functions.
Different cost functions relate to different driving objectives and consequently lead to different MFGs. In this subsection, we will present two concrete cost functions quantifying AVs' kinetic energy, driving efficiency and safety. 

\subsubsection{MFG-Separable}\label{sec:mfgsep}
We will propose a special cost function whose corresponding MFG has nice mathematical properties.
This family of cost functions is called \emph{separable} \cite{ambrose2018existence}, i.e., $f(u,\rho)$ can be written as the sum of two univariate functions with respect to $u$ and $\rho$.
Denote $\rho_{\text{jam}}$ the jam density, we propose a cost function that is separable and models AVs' kinetic energy, driving efficiency and safety:
\begin{align}
  f_{\text{Sep}}(u,\rho)=\underbrace{\frac12\left(\frac{u}{u_{\text{max}}}\right)^2}_{\text{kinetic energy}}-\underbrace{\frac{u}{u_{\text{max}}}}_{\text{efficiency}}+\underbrace{\frac{\rho}{\rho_{\text{jam}}}}_{\text{safety}}.\label{eq:costfct3}
\end{align}
In Eq.~\eqref{eq:costfct3}, the first term of $f_{\text{Sep}}(u,\rho)$ represents the kinetic energy; the second term quantifies driving efficiency by speed magnitudes; the last term quantifies driving safety using a traffic congestion penalty term on density $\rho$, meaning that AVs tend to avoid staying in high density areas. We denote the corresponding MFG by [MFG-Separable].

Since the cost function $f_{\text{Sep}}(u,\rho)$ is separable, [MFG-Separable] is a potential game when there are no speed constraints \cite{benamou2017variational}. \cite{cardaliaguet2015weak} proved the existence and uniqueness results for a family of potential MFGs including the one presented here.

When there are speed constraints $0\leq u\leq u_{\text{max}}$, the minimum of $f(u,\rho)+uV_x$ is attained at:
\begin{align}
  u=\max\left\{\min\left\{u_{\text{max}}(1-u_{\text{max}}V_x),u_{\text{max}}\right\},0\right\}.
\end{align}
So the MFG system is:
\begin{subequations}
\begin{numcases}{\text{[MFG-Separable]}\quad}
  \rho_t+(\rho u)_x=0,\label{eq:mfg_nocross_cont}\\
  V_t+uV_x+\frac12\left(\frac{u}{u_{\text{max}}}\right)^2-\frac{u}{u_{\text{max}}}+\frac{\rho}{\rho_{\text{jam}}}=0,\label{eq:mfg_nocross_hjb1}\\
  u=\max\left\{\min\left\{u_{\text{max}}(1-u_{\text{max}}V_x),u_{\text{max}}\right\},0\right\}.\label{eq:mfg_nocross_hjb2}
\end{numcases}
\end{subequations}

\subsubsection{MFG-NonSeparable}

We propose another cost function that quantifies driving safety in a more explicit way. The cost function is:
\begin{align}
  f_{\text{NonSep}}(u,\rho)=\underbrace{\frac12\left(\frac{u}{u_{\text{max}}}\right)^2}_{\text{kinetic energy}}-\underbrace{\frac{u}{u_{\text{max}}}}_{\text{efficiency}}+\underbrace{\frac{u\rho}{u_{\text{max}}\rho_{\text{jam}}}}_{\text{safety}}.\label{eq:costfct1}
\end{align}
It quantifies kinetic energy and driving efficiency in the same way as in $f_{\text{Sep}}$ but uses a different traffic congestion penalty term on the product of density and speed to quantify driving safety. The new penalty term means that AVs tend to decelerate in high density areas and accelerate in low density areas. We denote the corresponding MFG by [MFG-NonSeparable].

Let us rewrite the cost function in a different way:
\begin{align}
  f_{\text{NonSep}}(u,\rho)=\frac{1}{2u_{\text{max}}^2}(U(\rho)-u)^2-\frac12\left(1-\frac{\rho}{\rho_{\text{jam}}}\right)^2,\label{eq:costfct1n}
\end{align}
where $U(\rho)$ is the Greenshields desired speed function defined in Eq.~\eqref{eq:greenshields}. It provides a different way to interpret the cost function $f_{\text{NonSep}}$ that AVs tend to be not too far from human driving, and they like to stay in low density areas.

Comparing Eq.~\eqref{eq:costfct1n} with Eq.~\eqref{eq:costfun_lwr}, we see that $f_{\text{NonSep}}$ can be seen as a variant of $f_{\text{LWR}}$. That is one reason why we pick the cost function $f_{\text{LWR}}$ in Section~\ref{sec:lwrsolmfg}.

With speed constraints $0\leq u\leq u_{\text{max}}$, the minimum of $f(u,\rho)+uV_x$ is attained at:
\begin{align}
  u=\max\left\{\min\left\{u_{\text{max}}(1-\frac{\rho}{\rho_{\text{jam}}}-u_{\text{max}}V_x),u_{\text{max}}\right\},0\right\}.\label{eq:speednonsep}
\end{align}
The corresponding MFG system is:
{\small\begin{subequations}
\begin{numcases}{\text{[MFG-NonSeparable]}}
  \rho_t+(\rho u)_x=0,\label{eq:mfg_cross_cont}\\
  V_t+uV_x+\frac12\left(\frac{u}{u_{\text{max}}}\right)^2-\frac{u}{u_{\text{max}}}+\frac{u\rho}{u_{\text{max}}\rho_{\text{jam}}}=0,\label{eq:mfg_cross_hjb1}\\
  u=\max\left\{\min\left\{u_{\text{max}}(1-\frac{\rho}{\rho_{\text{jam}}}-u_{\text{max}}V_x),u_{\text{max}}\right\},0\right\}.\quad\quad\  \label{eq:mfg_cross_hjb2}
\end{numcases}
\end{subequations}}

\begin{remark}
  Letting $V_x\to0$, Eq.~\eqref{eq:speednonsep} becomes $u=u_{\text{max}}(1-\rho/\rho_{\text{jam}})$, which is the same as the Greenshields desired speed defined in Eq.~\eqref{eq:greenshields}. From Theorem~\ref{thm:myopic} we know that the Greenshields LWR is the myopic limit of [MFG-NonSeparable].
\end{remark}
{
[MFG-NonSeparable] can also be interpreted from the perspective of traditional higher-order traffic flow models. To see this, let $u_{\text{max}}=\rho_{\text{jam}}=1$, remove the speed constraints $0\leq u\leq u_{\text{max}}$ and take Eq.~\eqref{eq:mfg_cross_hjb2} into Eq.~\eqref{eq:mfg_cross_hjb1}. Then we obtain the following system:
\begin{subequations}
\begin{numcases}{}
  \rho_t+(\rho u)_x=0,\label{eq:appen1}\\
  V_t=\frac12u^2,\label{eq:appen2}\\
  \rho+u+V_x=1.\label{eq:appen3}
\end{numcases}
\end{subequations}
Differentiating Eq.~\eqref{eq:appen2} with respect to $x$ and Eq.~\eqref{eq:appen3} with respect to $t$ yields:
\begin{align}
  &V_{tx}=uu_x,\\
  &\rho_t+u_t+V_{xt}=0.
\end{align}
Using the identity $V_{xt}=V_{tx}$, we can eliminate the variable $V$ from the HJB equations and obtain:
\begin{align}
  u_t+uu_x=-\rho_t.\label{eq:appen4}
\end{align}
Eq.~\eqref{eq:appen4} coupled with the continuity equation \eqref{eq:appen1} forms the reduced MFG system:
\begin{subequations}
\begin{numcases}{\mbox{[reduced MFG]}\quad}
    \rho_t+(\rho u)_x=0,\label{eq:mfgac1}\\
    u_t+uu_x=-\rho_t.\label{eq:mfgac2}
\end{numcases}
\end{subequations}

The reduced MFG system has the same initial condition $\rho(x,0)=\rho_0(x)$ as that of the original system. Moreover, the original system's terminal condition $V_T(x)=C$ and Eq.~\eqref{eq:appen3} yield the terminal condition $\rho(x,T)+u(x,T)=1$ of [reduced MFG].

The reduced MFG system [reduced MFG] has a similar structure to traditional higher-order traffic flow models. Eq.~\eqref{eq:mfgac1} is the continuity equation and Eq.~\eqref{eq:mfgac2} has an interpretation that the car's acceleration $u_t+uu_x$ is exactly the negative temporal derivative of the density.
}

The proposed MFG systems are more effectively simulated than discrete differential games. We will discretize the systems in space and time and then present a solution algorithm to compute the MFE.

\section{MFE solution algorithm}\label{sec:sol}

Because of the forward-backward structure, the MFG system can be solved in neither forward nor backward direction. Given the density profile $\rho(x,t)$, the HJB equations~\eqref{eq:mfg_hjb1}\eqref{eq:mfg_hjb2} can be solved backward from $t=T$ to $t=0$ with terminal cost $V_T(x)$ for $u(x,t)$ and $V(x,t)$; given the velocity field $u(x,t)$, the continuity equation~\eqref{mfg_cont} can be solved forward from $t=0$ to $t=T$ with initial density $\rho_0(x)$ for $\rho(x,t)$. However, the two directions can not be matched simultaneously. So it is challenging to compute the MFE numerically.

Based on the existing studies, there have been three types of numerical methods  considered for MFG: fixed-point iteration, variational method and Newton's method.

The fixed-point iteration solves the forward and backward equations alternatingly. It is easy to implement once appropriate forward and backward solvers are picked \cite{couillet2012electrical,chevalier2015micro}.
However, the iterations converge only when $T$ is small, that is, for a short planning horizon. Moreover, there is no theory to estimate how small $T$ should be to guarantee the convergence.

The variational method deals with separable cost functions and potential MFGs. In this case, it is shown that the MFG system is equivalent to an optimization problem constrained by the continuity equation \cite{benamou2017variational}. Then a variety of optimization tools can be applied \cite{lachapelle2011mean,benamou2015augmented,chow2018algorithm}. The variational method works for any planning horizon but relies on the separability of the cost function. \cite{lachapelle2011mean} used the variational method to solve MFGs in pedestrian crowds modeling.

A more general approach is based on the Newton's method. Such an approach is first proposed by \cite{achdou2010mean,achdou2012mean,achdou2012iterative} to solve a family of MFGs. The key idea is to take both forward and backward equations as a single nonlinear system and solve the nonlinear system by the Newton's method. This method is suitable for our purpose because it has no requirements on the length of the planning horizon nor the separability of the cost function. However, the Newton's method may fail to converge if one does not have a good initial guess to the solution. So tricks to improve the convergence are needed when applying the Newton's method.

This paper develops a multigrid preconditioned Newton's finite difference algorithm for MFG. It works well with different cost functions and planning horizons. Numerical examples of MFGs proposed in Section \ref{sec:flow} are shown using this algorithm.


\subsection{Algorithm}
Let us divide the road $[0,L]$ into cells $\{[x_{j-1},x_j]\}_{j=1}^{N_x}$ and the planning horizon $[0,T]$ into time steps $\{t^n\}_{n=0}^{N_t}$ with spatial and temporal step sizes $\Delta x=L/N_x$ and $\Delta t=T/N_t$. To impose the periodic boundary conditions, $x_0$ and $x_{N_x}$ are assumed to be the same location. Denote $\rho_j^n$ the average density and $u_j^n$ the average velocity on the $j^{\text{th}}$ cell $[x_{j-1},x_j]$ at time $t^n$ ($j=1,\dots,N_x$; $n=0,\dots,N_t$). Denote $V_j^n=V(x_j,t^n)$.

We first discretize the continuity equation~\eqref{mfg_cont} by a finite volume conservative Lax-Friedrichs scheme \cite{leveque2002finite}:
\begin{align}
  \rho_{j}^{n+1}=\frac12(\rho_{j-1}^n+\rho_{j+1}^n)-\frac{\Delta t}{2\Delta x}(\rho_{j+1}^nu_{j+1}^n-\rho_{j-1}^nu_{j-1}^n).\label{eq:schemecont}
\end{align}

We then discretize the HJB equations~\eqref{eq:mfg_hjb1}\eqref{eq:mfg_hjb2} by an upwind scheme:
\begin{align}
  &\frac{V_j^{n+1}-V_j^n}{\Delta t}+f^*\left(\frac{V_{j}^{n+1}-V_{j-1}^{n+1}}{\Delta x},\rho_j^n\right)=0,\label{eq:schemehjb1}\\
  &u_j^n=f^*_p\left(\frac{V_{j}^{n+1}-V_{j-1}^{n+1}}{\Delta x},\rho_j^n\right).\label{eq:schemehjb2}
\end{align}

\begin{remark}
  To ensure the stability of scheme~\eqref{eq:schemecont}, the CFL condition $\alpha\Delta t\leq\Delta x$ should be posed \cite{leveque2002finite} where $\alpha=\max\nolimits_{j,n}|u_j^n|$. When the MFG has speed constraints $0\leq u\leq u_{\text{max}}$, it suffices to ensure $u_{\text{max}}\Delta t\leq\Delta x$.
\end{remark}

The initial and terminal conditions are discretized by:
\begin{align}
  \rho_j^0=\frac{1}{\Delta x}\int_{x_{j-1}}^{x_j} \rho_0(x)\,dx,\quad V_j^{N_t}=V_T(x_j).\ \ (j=1,\dots,N_x)\label{eq:schemeinitercond}
\end{align}

Eqs.~\eqref{eq:schemecont}\eqref{eq:schemehjb1}\eqref{eq:schemehjb2}\eqref{eq:schemeinitercond} form a closed system for unknowns $\{\rho_j^n\}_{1\leq j\leq N_x}^{0\leq n\leq N_t}$,\break $\{u_j^n\}_{1\leq j\leq N_x}^{0\leq n\leq N_t-1}$ and $\{V_j^n\}_{1\leq j\leq N_x}^{0\leq n\leq N_t}$. The system can be written as:
\begin{align}
  F(w)=0\label{eq:nonlinsys},
\end{align}
where $w\in\mathbb{R}^{3N_xN_t+2N_x}$ is a long vector containing all $\rho_j^n$, $u_j^n$ and $V_j^n$, and $F:\mathbb{R}^{3N_xN_t+2N_x}\to\mathbb{R}^{3N_xN_t+2N_x}$ encodes all equations.

Eq.~\eqref{eq:nonlinsys} may lead to a large nonlinear system. We denote $J$ the Jacobian matrix of $F$ and apply the Newton's method to solve Eq.~\eqref{eq:nonlinsys}:
\begin{align}
  w^{n+1}=w^n-J(w^n)^{-1}F(w^n),
\end{align}
with any initial guess $w^0$.

To improve the convergence of Newton's iterations, we apply multigrid to get a good initial guess and preconditioning to accelerate the linear solver. Multigrid and preconditioning are widely used tricks in numerical algorithms, see \cite{hackbusch2013multi,golub2012matrix}. 
\begin{itemize}
  \item \textbf{Multigrid}: Start with a coarse grid $N_x^{(0)}\times N_t^{(0)}$ so that the MFG system is easy to solve. Then iteratively refine the grids and solve the MFG system on finer grids $N_x^{(k)}\times N_t^{(k)}$, $k=1,2,\dots$ until getting a solution of desired resolution. At step $k$, interpolate the solution $w^{(k-1)}$ from the grid $N_x^{(k-1)}\times N_t^{(k-1)}$ onto the finer grid $N_x^{(k)}\times N_t^{(k)}$, which provides a good initial guess when solving on the finer grid by the Newton's method.
  \item \textbf{Preconditioning}: At each Newton's iteration a linear system
  \begin{align}
    J(w^n)(w^{n+1}-w^n)=-F(w^n),
  \end{align} need to be solved. We use the GMRES iterative linear solver \cite{saad1986gmres} since $J(w^n)$ is sparse . However, the ill-posedness of the linear system leads to bad convergence. To solve the issue, we pick an approximate matrix $\tilde{J}(w^n)$ to $J(w^n)$ by ignoring the coupling parts between forward and backward equations. Inverting $\tilde{J}(w^n)$ is equivalent to solving a decoupled forward-backward system. We use $\tilde{J}(w^n)$ as a preconditioner to accelerate the GMRES convergence.
\end{itemize}

Using the algorithm, we shall compute MFE solutions and show simulations for MFGs proposed in Section~\ref{sec:flow}.

\subsection{Numerical examples}

\subsubsection{Settings}

Set the road length $L=1$ and the planning horizon length $T=3$. Set the free flow speed $u_{\text{max}}=1$ and the jam density $\rho_{\text{jam}}=1$. Choose the following initial density:
\begin{align}
  \rho_0(x)=\rho_a+(\rho_b-\rho_a)\exp[-\frac{(x-L/2)^2}{2\gamma^2}], \label{eq:simuini}
\end{align}
where $0\leq\rho_a\leq\rho_b\leq1$ and $\gamma>0$ are constant parameters.
The initial density represents the scenario that initially cars cluster near $x=L/2$ and the traffic is lighter in other places. We choose the terminal cost $V_T(x)=0$ and specify periodic boundary conditions $\rho(0,t)=\rho(L,t)$, $V(0,t)=V(L,t)$. 


 On a spatial-temporal grid of size $N_x=120$ and $N_t=480$, we compute MFE solutions $\rho^*(x,t)$, $u^*(x,t)$ and $V^*(x,t)$ for the three MFG systems in Section~\ref{sec:flow}.

For [MFG-LWR], we choose $U(\rho)$ to be the Greenshields desired speed function defined in Eq.~\eqref{eq:greenshields}.
\subsubsection{Density Evolution}
\begin{figure}[htbp]
  \centering
  \subfloat[\mbox{[MFG-NonSeparable]}]{\label{fig:densityevolutioncross}
    \includegraphics[width=.48\textwidth]{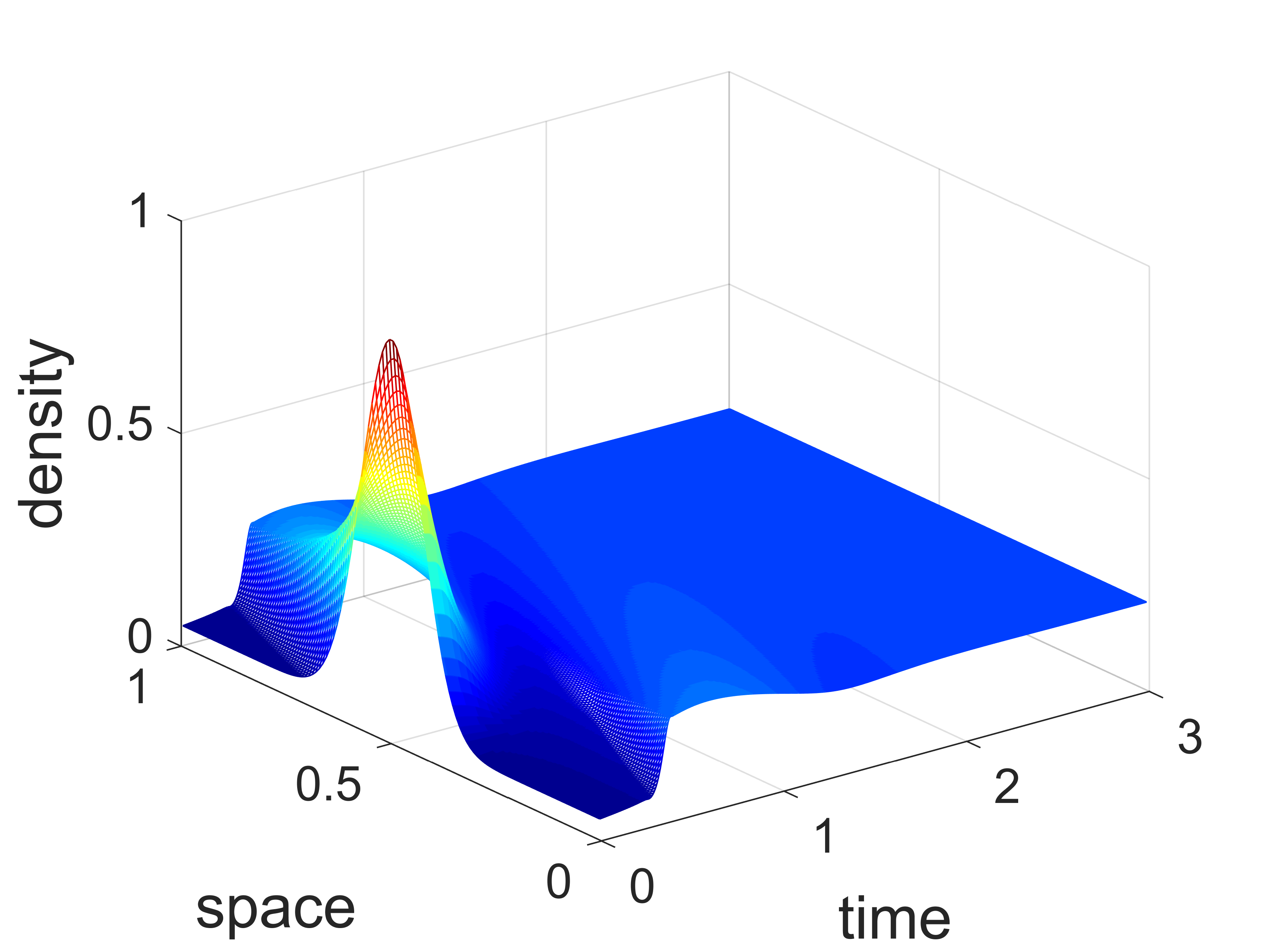}
  }
  \subfloat[\mbox{[MFG-Separable]}]{
    \includegraphics[width=.48\textwidth]{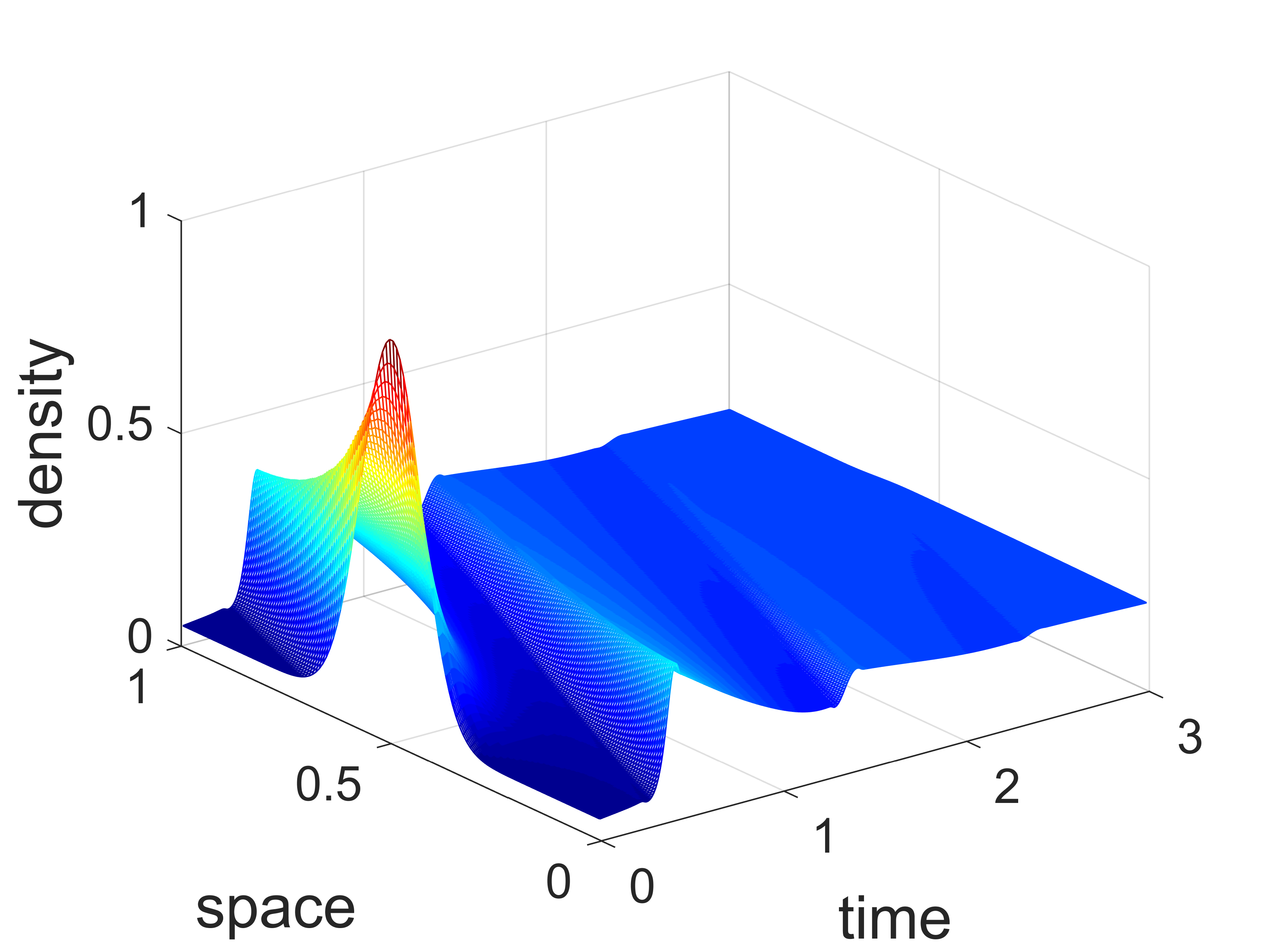}
  }
  \caption{Density Evolution of [MFG-NonSeparable] and [MFG-Separable]}
  \label{fig:densityevolution}
\end{figure}

Fix the same initial density $\rho_0(x)$ defined in Eq.~\eqref{eq:simuini} with $\rho_a=0.05$, $\rho_b=0.95$ and $\gamma=0.1$, we compute the MFE solutions for the three MFG systems and plot their density evolutions in a 3D space-time-density diagram. See Figure~\ref{fig:mfg-lwra} and Figure~\ref{fig:densityevolution}.

Figure~\ref{fig:mfg-lwra} shows the formation and propagation of a shock wave for [MFG-LWR]. The shock wave moves with smaller and smaller amplitudes but does not disappear in the given time horizon $[0,T]$.

Figure~\ref{fig:densityevolution} shows that for both [MFG-NonSeparable] and [MFG-Separable], the initial high density quickly dissipates. For [MFG-NonSeparable], the density profile keeps smooth and no shock wave forms. From time $t=1$, the density becomes a uniform flow. For [MFG-Separable], the behavior is similar but the high density dissipates more slowly and the density becomes nearly a uniform flow from $t=2$. Such phenomena are different from traditional traffic flow models.

\begin{figure}[htbp]
  \centering
  \subfloat[\mbox{[MFG-NonSeparable]} at $t=0$]{\label{fig:rhouV_t0_cross}
    \includegraphics[width=.48\textwidth]{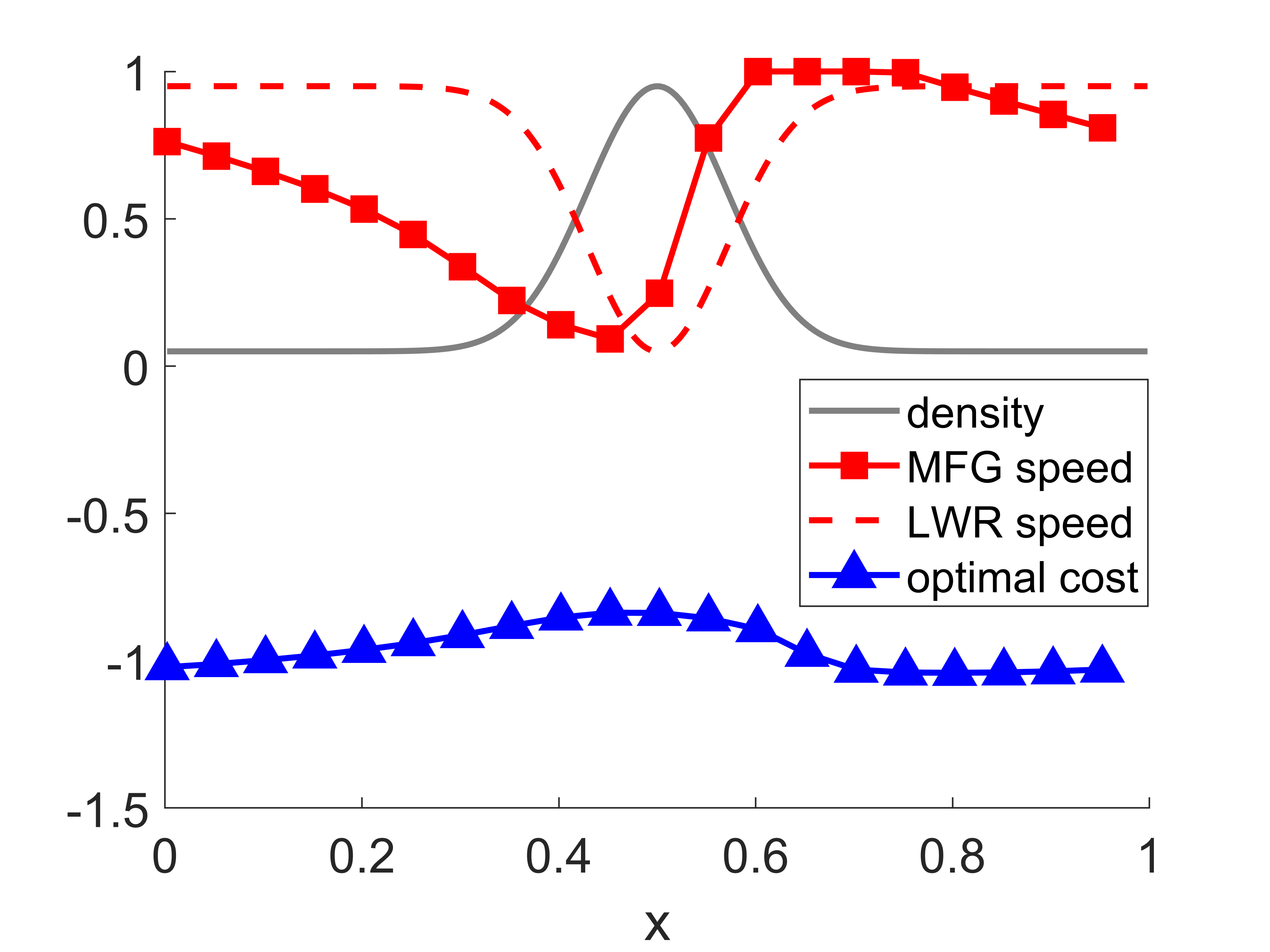}
  }
  \subfloat[\mbox{[MFG-Separable]} at $t=0$]{\label{fig:rhouV_t0_nocross}
    \includegraphics[width=.48\textwidth]{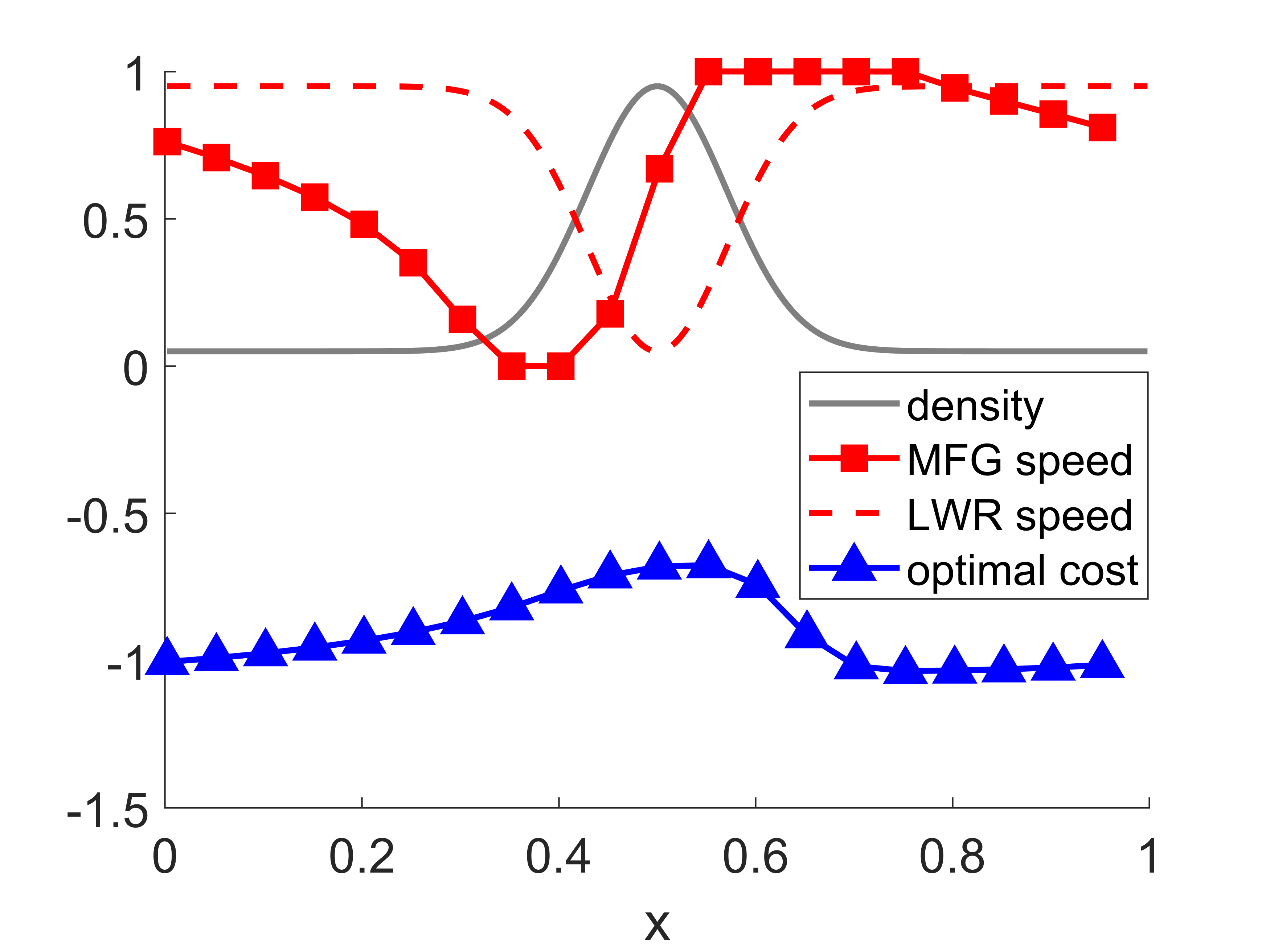}
  }

  \subfloat[\mbox{[MFG-NonSeparable]} at $t=1.5$]{\label{fig:rhouV_t1_cross}
    \includegraphics[width=.48\textwidth]{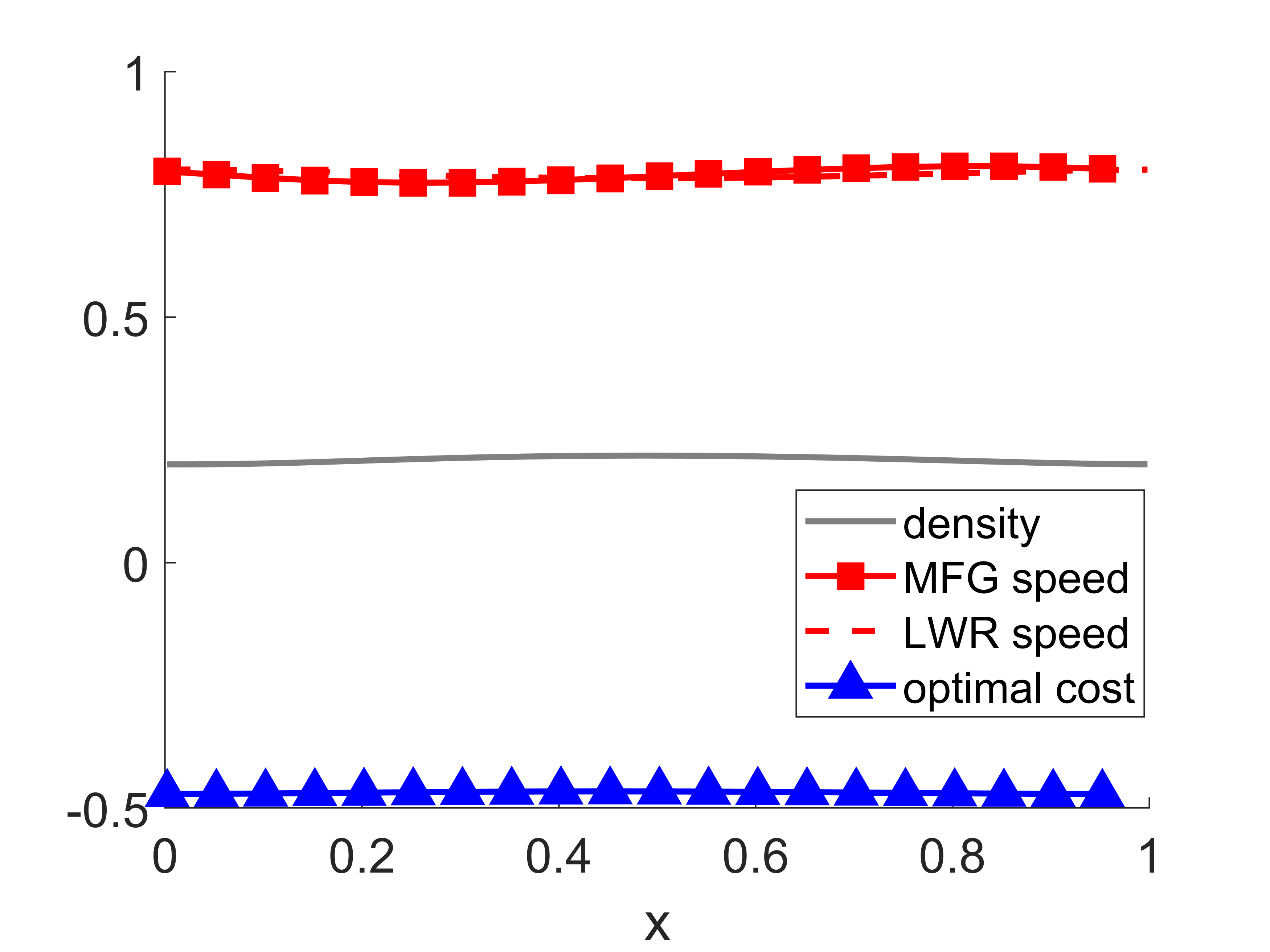}
  }
  \subfloat[\mbox{[MFG-Separable]} at $t=1.5$]{\label{fig:rhouV_t1_nocross}
    \includegraphics[width=.48\textwidth]{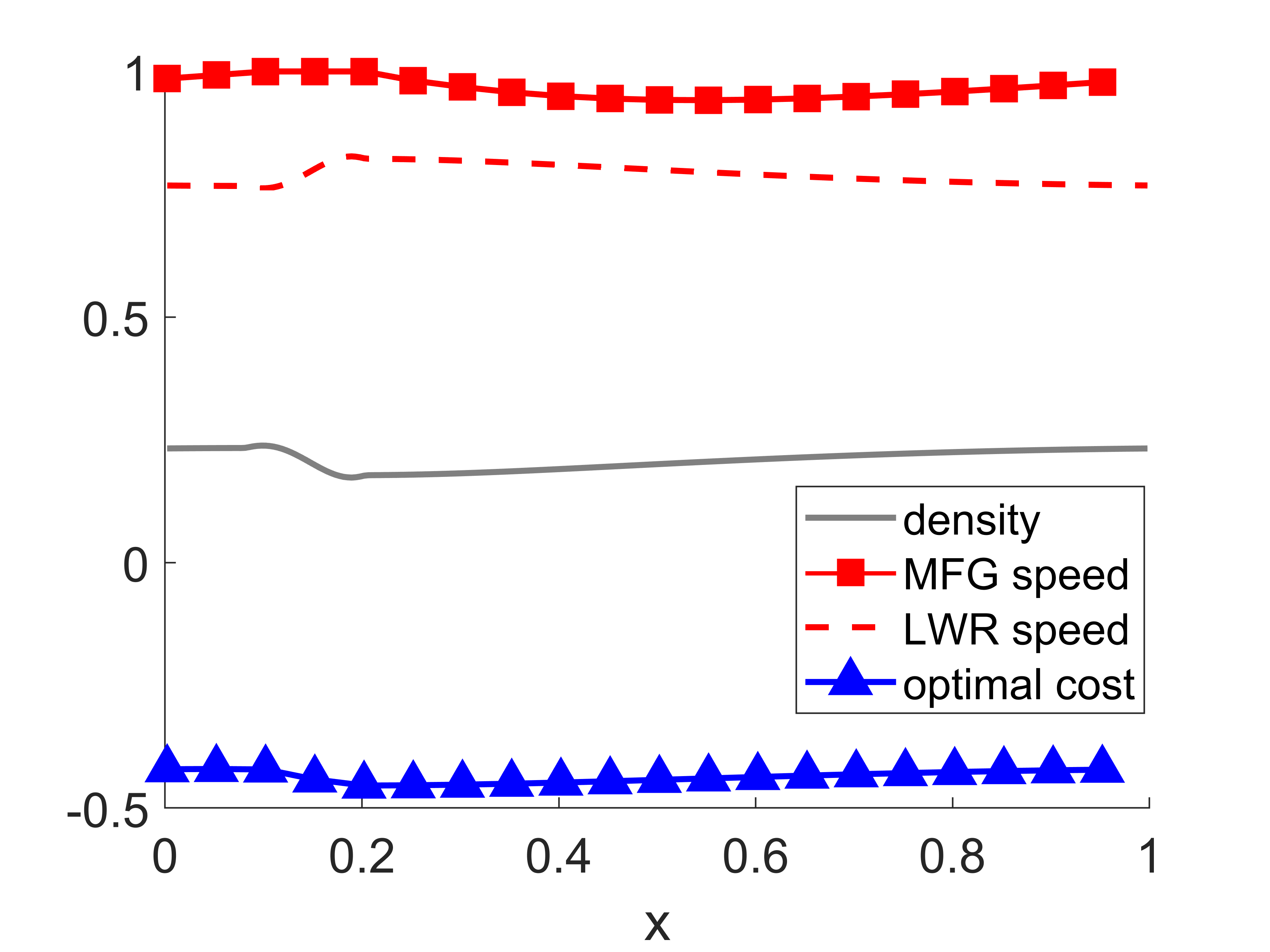}
  }

  \caption{Density, speed and optimal cost profiles for [MFG-NonSeparable] and [MFG-Separable] at $t=0$ and $t=1.5$}
  \label{fig:rhouV}
\end{figure}

The results show that AVs' anticipation behavior helps to avoid the formation of shock waves and to stabilize the traffic in this set-up.
Figure~\ref{fig:rhouV} reveals the rationale by plotting the snapshots of the MFE solutions' density, speed and optimal cost profiles with respect to spatial coordinate $x$ at time instants $t=0$ and $t=1.5$ for [MFG-NonSeparable] and [MFG-Separable]. In addition, we compute the LWR speeds from the density profiles with the Greenshields desired speed function defined in Eq.~\eqref{eq:greenshields}, plot the LWR speeds and compare them with the MFE speeds in the same axes.

We observe from Figure~\ref{fig:rhouV} that all of the density, speed and optimal cost profiles converge to constant profiles as time goes on. Figure~\ref{fig:rhouV_t0_cross} and Figure~\ref{fig:rhouV_t0_nocross} show asymmetric optimal cost around a jam area with symmetric traffic density. As a result of the ``pressure'' from the optimal cost,
AVs tend to slow down farther upstream before joining the jam and immediately speed up after leaving the jam,
in contrast to HVs whose speeds remain symmetric before and after the jam area.
In other words,
a HV's speed is determined only through traffic density at the current time 
while that of an AV depends on model predicted traffic density over the entire horizon. Such behavioral difference between AVs and HVs result in different traffic flows.

\subsubsection{Fundamental Diagram}
Fundamental diagram is a basic tool to understand traditional traffic flow models \cite{orosz2010traffic}. In this subsection we will collect density and flow data from the MFE solution and plot the data in a fundamental diagram for both [MFG-LWR] and [MFG-NonSeparable].

The density and flow data are collected from the MFE solution as follows: Take $n_x=24$ equidistantly
distributed locations $x_1,x_2,\dots,x_{n_x}$ on $[0,L]$ and $n_t=96$ time snapshots $t^k=\frac{kT}{n_t}$, $k=0,1,\dots,n_t$. We first pick the density and speed values $\rho^*(x_i,t^k)$ and $u^*(x_i,t^k)$ that represent the average density and speed near the spatial-temporal coordinate $(x_i,t^k)$ for $i=1,\dots,n_x$, $k=0,\dots,n_t$, then compute the flow $q^*(x_i,t^k)=\rho^*(x_i,t^k)u^*(x_i,t^k)$ from the density and speed. The collected data $\{\rho^*(x_i,t^k),q^*(x_i,t^k)\}_{1\leq i\leq n_x,0\leq k\leq n_t}$ are plotted on the density-flow diagram. Such a way to plot the fundamental diagram from a macroscopic traffic flow model is also used in \cite{siebel2006fundamental}.

\begin{figure}[htbp]
  \centering
  \subfloat[Density Evolution]{\label{fig:mfg-lwra}
    \includegraphics[width=.48\textwidth]{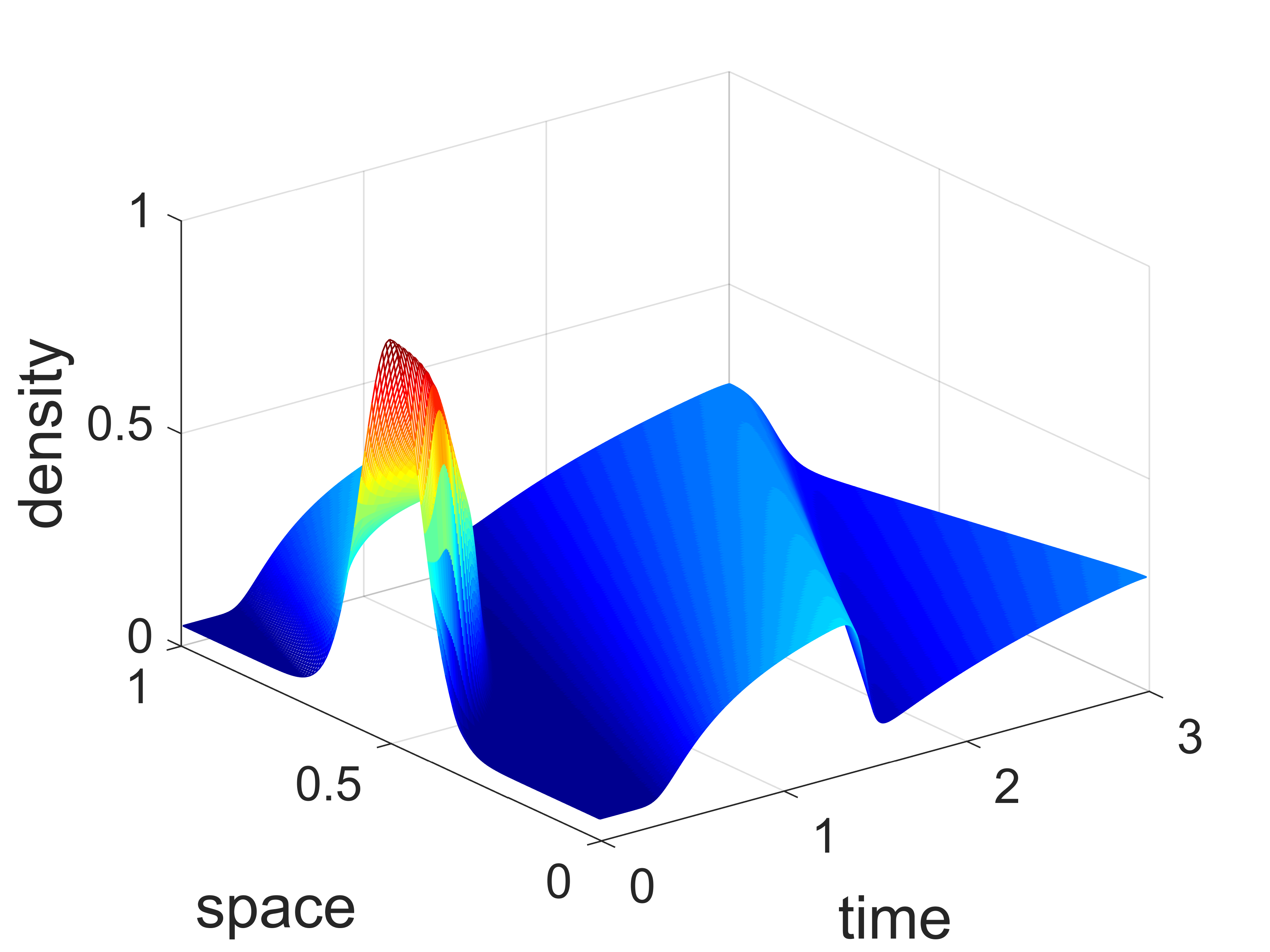}
  }
  \subfloat[Fundamental Diagram]{\label{fig:mfg-lwrb}
    \includegraphics[width=.48\textwidth]{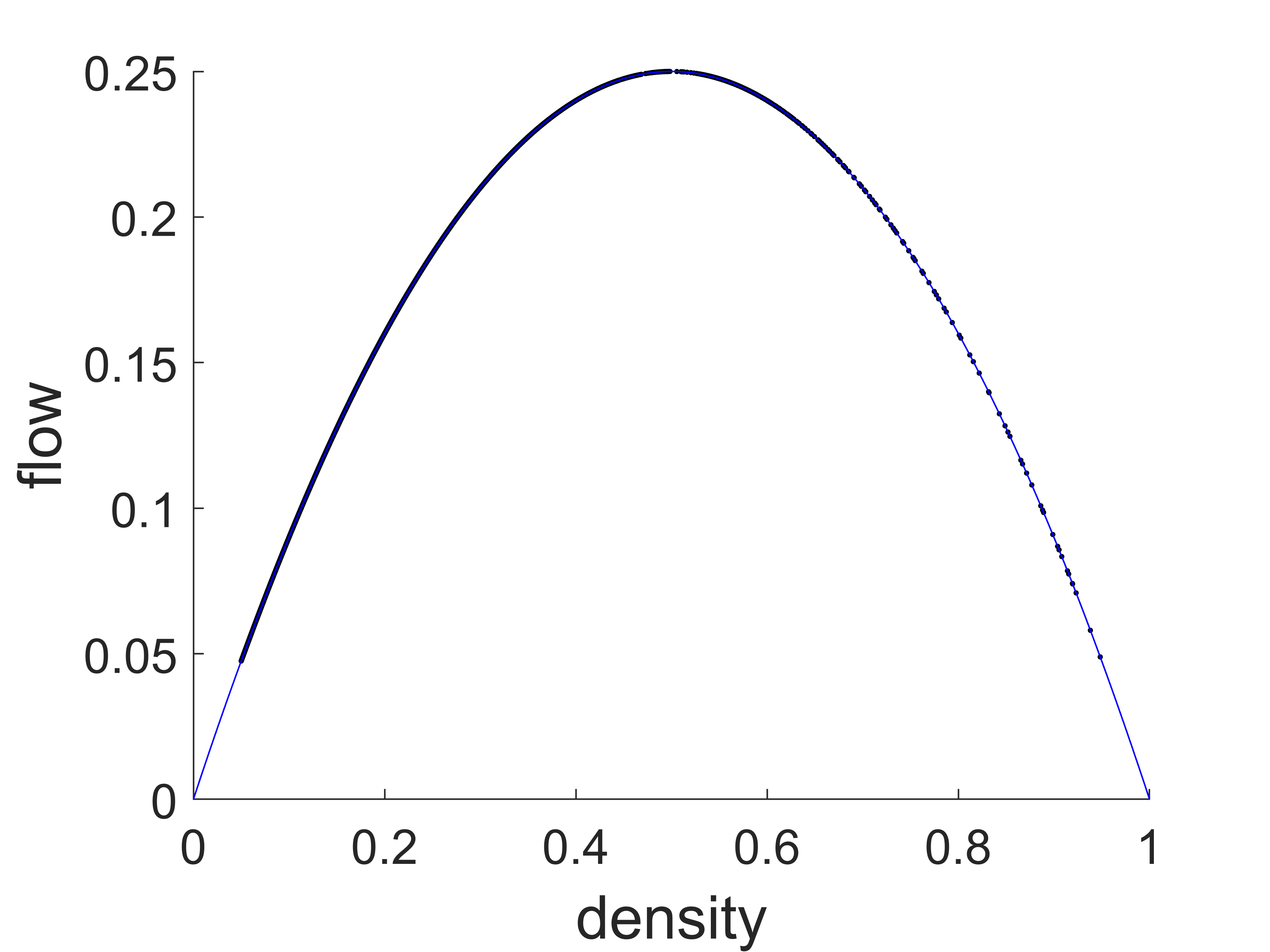}
  }
  \caption{[MFG-LWR]}
  \label{fig:mfg-lwr}
\end{figure}

For [MFG-LWR], we collect the data from the MFE solution shown in Figure~\ref{fig:mfg-lwra} and plot the fundamental diagram in Figure~\ref{fig:mfg-lwrb}. We see that all of the collected density-flow data points fall onto the Greenshields equilibrium curve $q=u_{\text{max}}\rho(1-\rho/\rho_{\text{jam}})$.
The results verify Theorem~\ref{thm:mfglwr}.

\begin{figure}[htbp]
  \centering
  \includegraphics[width=.75\textwidth]{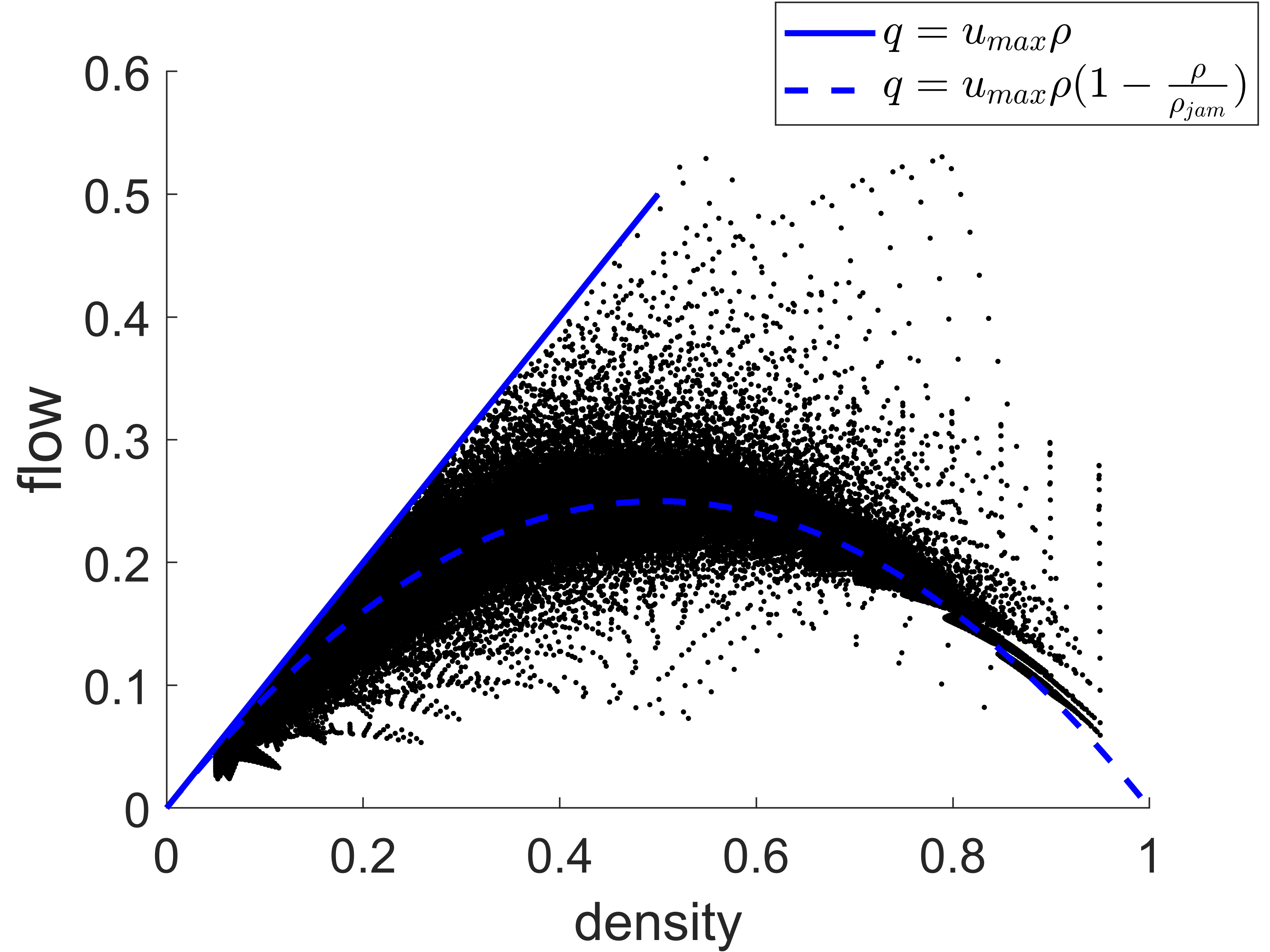}
  \caption{Fundamental diagram of [MFG-NonSeparable]}
  \label{fig:fdmfg}
\end{figure}

For [MFG-NonSeparable], we plot the fundamental diagram by collecting the data from a set of different MFE solutions. We take different initial densities by varying the values of $\rho_a$ and $\rho_b$ from $0.05$ to $0.95$ but keep $\rho_a<\rho_b$ and $\gamma=0.1$. For each initial density $\rho_0(x)$ we compute the MFE solution from [MFG-NonSeparable] and collect the data in the way mentioned. Then we plot all collected data in the same fundamental diagram, see Figure~\ref{fig:fdmfg}.

We observe from Figure~\ref{fig:fdmfg} that: (i) All data points lie below the line $q=u_{\text{max}}\rho$, this is due to the speed constraint $u\leq u_{\text{max}}$ in [MFG-NonSeparable]. (ii) All data points cluster around the Greenshields equilibrium curve $q=u_{\text{max}}\rho(1-\rho/\rho_{\text{jam}})$, this is because [MFG-NonSeparable] is related to the Greenshields LWR. The traffic flow always converges to a uniform flow represented by a data point on the Greenshields equilibrium curve. The position of the data point on the curve depends on the initial density. (iii) The fundamental diagram can be split into a free flow regime where $\rho\leq0.5$ and a congested regime where $\rho>0.5$. In the free flow regime, cars can achieve the free flow speed $u_{\text{max}}$. It reflects the efficiency term in the cost function $f_{\text{NonSep}}$; in the congested regime cars cannot achieve the free flow speed, which reflects the safety term in the cost function.  (iv) Different from the fundamental diagram of human driving \cite{seibold2012constructing}, data points in the free flow regime of Figure~\ref{fig:fdmfg} may not lie on the line $q=u_{\text{max}}\rho$, this results from AVs' anticipation behavior. Even in a low density area, the car may drive at a lower speed than the desired one if there is a traffic jam ahead.


\begin{figure}
  \centering
  \includegraphics[width=.5\textwidth]{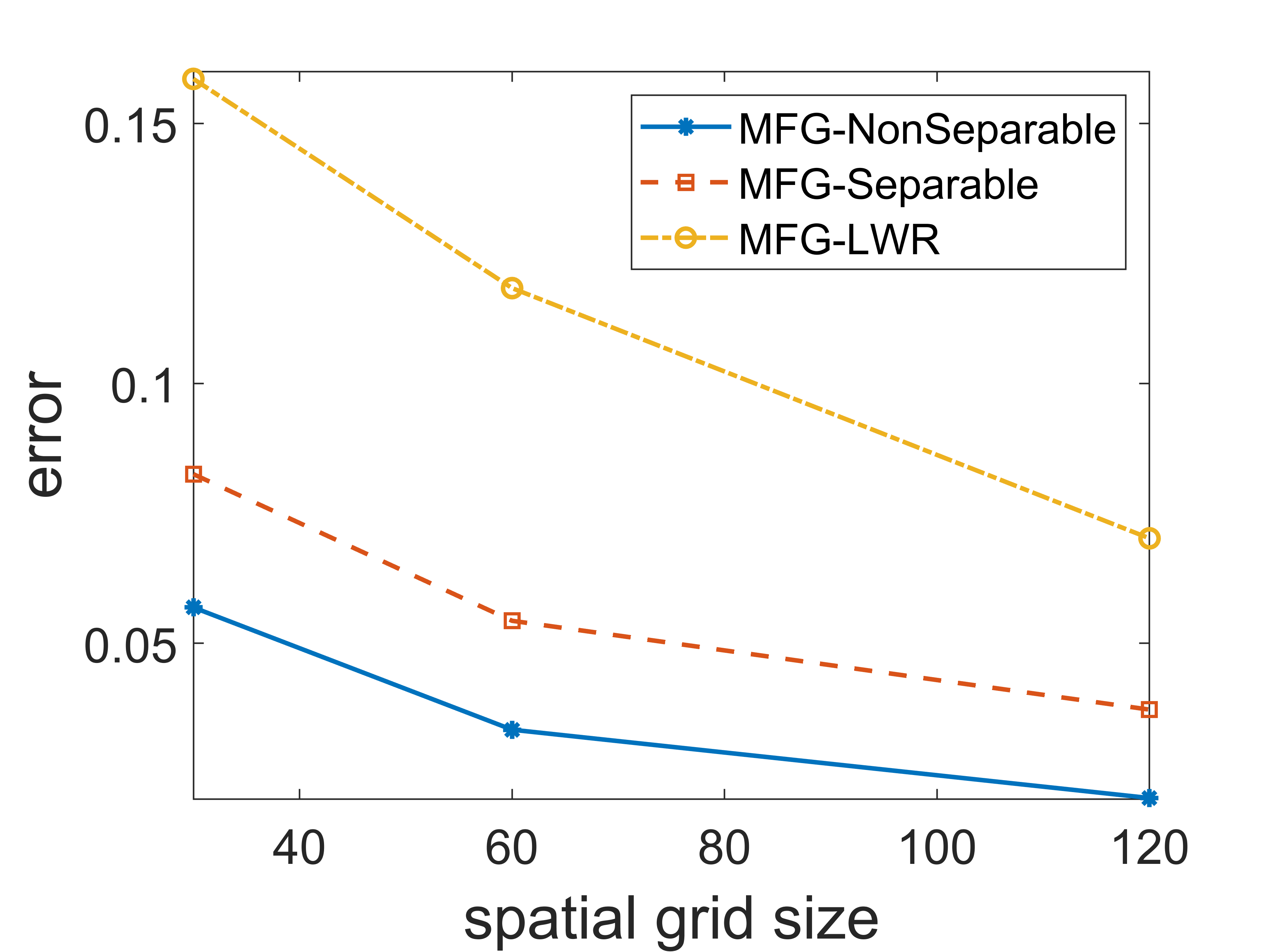}
  \caption{Convergence of solution algorithm in $L^1$ norm}
  \label{fig:algconvergence}
\end{figure}

\subsubsection{Algorithm Convergence}
To provide more evidence on the convergence of the solution algorithm, we compute and plot the solution errors on different grids for above examples. See Figure~\ref{fig:algconvergence}. Since we do not know any explicit solutions to those MFG systems, the errors are estimated using numerical solutions on different grids. To check the solution $\rho^*(x,t),u^*(x,t)$ on the $N_x\times N_t$ grid, we first solve the MFG on the coarse grid of size $N_x/2\times N_t/2$ and then interpolate the coarse solution back onto the $N_x\times N_t$ grid. Denote the interpolated solution by $(\tilde{\rho}^*,\tilde{u}^*)$, the solution error on the $N_x\times N_t$ grid is estimated as:
\begin{align}
  \text{Error}_{N_x\times N_t}=\norm{\rho^*-\tilde{\rho}^*}+\norm{u^*-\tilde{u}^*},\label{eq:converr}
\end{align}
where the norm is chosen as the $L^1$ norm on $[0,L]\times[0,T]$.

We fix the spatial-temporal ratio $N_t/N_x=4$ and increase $N_x$ from 30 to 120. Then we plot the errors computed by Eq.~\eqref{eq:converr} with the spatial grid size $N_x$. From Figure~\ref{fig:algconvergence} we see first order convergence for all of the three numerical examples.

\section{From MFG back to \textit{N}-Car differential game}\label{sec:micro-macro}

Summarizing the previous sections, we have derived a continuous mean field game [MFG] from a discrete differential game (Definition~\ref{def:dg}) and developed a solution algorithm for the mean field game.
In this section we shall build the connection between the discrete differential game equilibrium (DGE) and the continuous mean field equilibrium (MFE) in the sense of $\epsilon$-Nash equilibrium. First we provide a way to construct a tuple of discrete controls from a MFE solution. Then we introduce the concept of $\epsilon$-Nash equilibrium and show how to characterize the accuracy of the MFE-constructed controls. It is validated by numerical examples that the MFE-constructed controls are a good approximate equilibrium of the original $N$-car differential game (DG) when $N$ is large.


\subsection{MFE-constructed controls and accuracy characterization}\label{sec:construction}

From a continuous MFE solution, we construct a tuple of discrete controls $\hat{v}_1(t),\dots,\hat{v}_N(t)$ for the DG by
applying the feedback law \eqref{eq:feed1}\eqref{eq:feed2}.
The rationale underlying such construction is quite straightforward: for $i=1,\cdots, N$,
the $i^{\text{th}}$ car's instantaneous speed selection at time $t$ is determined by MFE's velocity field $u^{\ast}(x_i(t),t)$ at that time and the $i^{\text{th}}$ car's location $x_i(t)$.
Mathematically, for $i=1,\cdots, N$ and $t\in \left[0, T\right]$:
\begin{align}
&\hat{v}_i(t)=u^{\ast}(x_i(t),t), \label{eq:dyn2}\\
&\dot{x}_i(t)=\hat{v}_i(t),\quad x_i(0)=x_{i,0}. \label{eq:dyn1}
\end{align}

Integrating the above dynamical system gives the $i^{\text{th}}$ car's velocity control $\hat{v}_i(t)$ and trajectory $x_i(t)$ over the planning horizon $[0, T]$, $i=1,\dots,N$.

\begin{figure}[htbp]
  \centering
  \includegraphics[width=.5\textwidth]{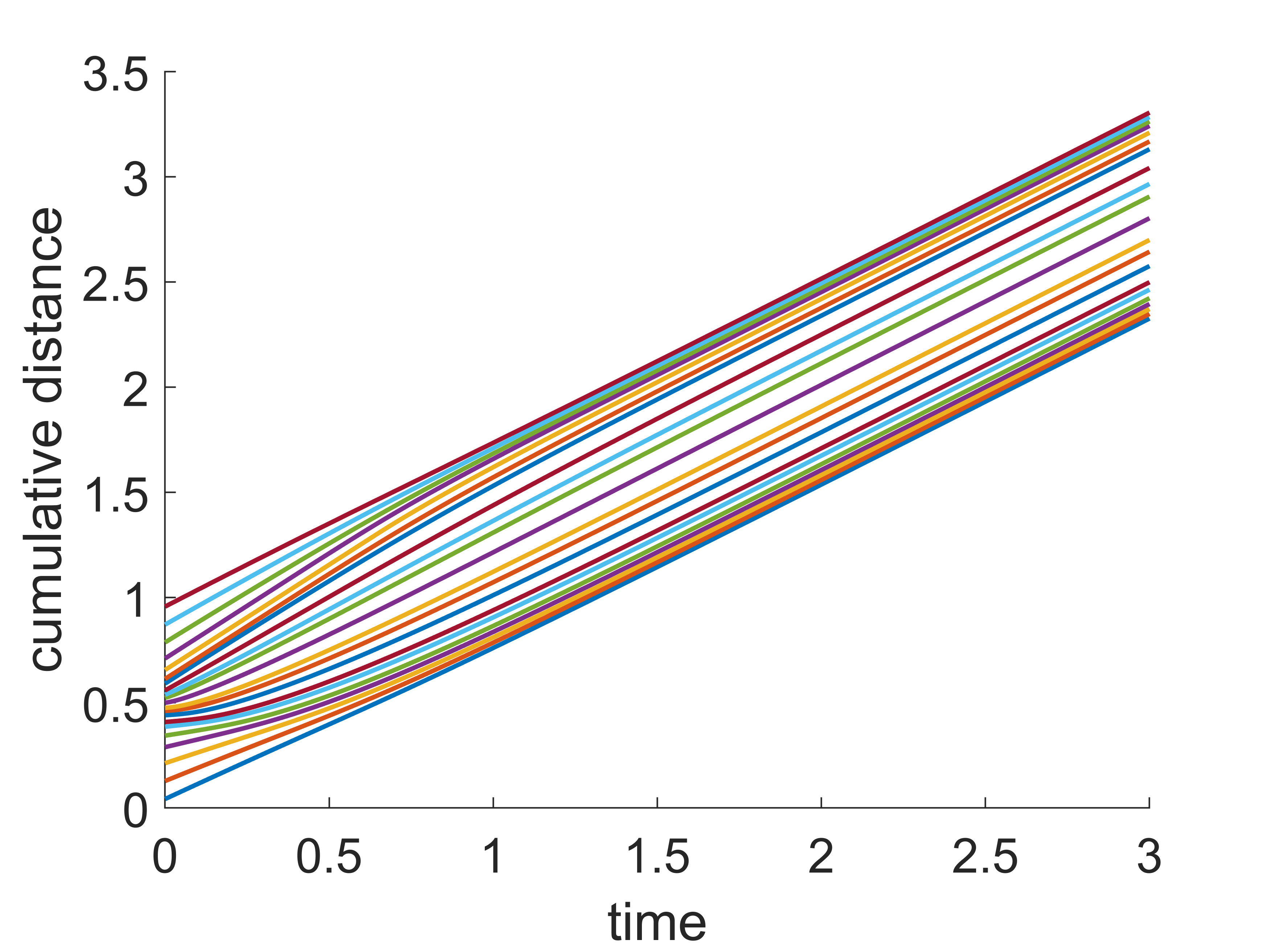}
  \caption{$N=21$ cars' trajectories integrated from the MFE solution of [MFG-NonSeparable]}
  \label{fig:traj}
\end{figure}

As an example, Figure~\ref{fig:traj} shows the trajectories integrated from the MFE solution shown in Figure~\ref{fig:densityevolutioncross} of [MFG-NonSeparable] with $N=21$ cars. Each of the 21 lines represents the trajectory of one car. We take the time to be the $x$-axis and the cumulative distance to be the $y$-axis to avoid special treatments on periodic boundary conditions. We see that even though cars cluster near $x=0.5$ at the staring time, they become uniformly distributed at the final time. It means that the flow converges to the uniform flow.

\begin{remark}
  We observe from Figure~\ref{fig:traj} that there is no intersection between any pair from the cars' trajectories. In other words, the first-in-first-out (FIFO) property is satisfied. Actually, Eqs.~\eqref{eq:dyn2}\eqref{eq:dyn1} is a first order ODE system for all $i=1,\dots,N$. The solution uniqueness then guarantees that there is no intersection between any pair from the $N$ trajectories starting from different initial locations \cite{strogatz2018nonlinear}. In other words, the FIFO property is always guaranteed by MFE-constructed controls.
\end{remark}

Now we would like to know whether the MFE-constructed controls $\hat{v}_1(t),\dots,\hat{v}_N(t)$ are a good approximate equilibrium of the DG. Since the DG's true equilibrium $v^{\ast}_1(t),\dots,v^{\ast}_N(t)$ may not exist nor be unique and typically it is hard to get,
we will characterize the accuracy of the constructed controls in terms of the driving cost functional.
Along this line, such an approximate equilibrium can be treated as an $\epsilon$-Nash equilibrium of the DG, which is formally defined below \cite{cardaliaguet2010notes,cardaliaguet2015weak}.
\begin{definition}
  A tuple of controls $\tilde{v}_1,\dots,\tilde{v}_N$ is an $\epsilon$-Nash equilibrium of the DG, if
  \begin{align}\label{def:epsNE}
    J_i^N(\tilde{v}_i,\tilde{\bm{v}}_{-i}) \leq J_i^N(v_i,\tilde{\bm{v}}_{-i}) + \epsilon,\quad \forall v_i\in\mathcal{A}, \quad i=1,\dots,N.
  \end{align}
\end{definition}
At an $\epsilon$-Nash equilibrium, no car can improve its driving cost better than $\epsilon$ by unilaterally switching its velocity control.

For a potential game where the cost function $f(u,\rho)$ is separable, e.g., [MFG-Separable], 
\cite{cardaliaguet2015weak} proved the correspondence between the MFE-constructed controls and an $\epsilon$-Nash equilibrium of DG, that is:
for any $\epsilon>0$, there exist $N, \sigma>0$ such that $\hat{v}_1(t),\dots,\hat{v}_N(t)$ constructed from the MFE solution 
is an $\epsilon$-Nash equilibrium of the DG.


Unfortunately, there do not exist any theoretical results for a general cost function such as $f_{\text{NonSep}}$.
In this paper, instead of offering a formal proof, we will validate such correspondence using numerical examples.
A rigorous proof of such correspondence for a general cost function will be left for future research.

Tailoring to our context, in the subsequent numerical examples, we aim to illustrate that the MFE-constructed controls are an $\epsilon$-Nash equilibrium of DG
by characterizing the accuracy $\epsilon$ across all feasible controls and all cars.
There always exists an arbitrarily large $\epsilon$ that can make $\hat{v}_1(t),\dots,\hat{v}_N(t)$ satisfy the condition \eqref{def:epsNE}.
What we are more interested in is a lower bound, denoted as $\hat{\epsilon} \geq 0$, such that:
\begin{equation}
   J_i^N(\hat{v}_i, \hat{\bm{v}}_{-i}) \leq  J_i^N(v_i, \hat{\bm{v}}_{-i}) + \hat{\epsilon},\quad \forall v_i\in\mathcal{A},\quad i=1,\dots,N.
\end{equation}


In other words,
\begin{equation}
  \hat{\epsilon} \triangleq \max_{1\leq i\leq N} \max_{v_i\in\mathcal{A}} \left\lbrace J_i^N(\hat{v}_i, \hat{\bm{v}}_{-i}) - J_i^N(v_i, \hat{\bm{v}}_{-i}) \right\rbrace.  \label{eq:fidelity}
\end{equation}
Let us move the second maximum symbol in front of the second term that depends only on $v_i$, then we have:
\begin{equation}
  \hat{\epsilon} = \max_{1\leq i\leq N} \left\lbrace J_i^N(\hat{v}_i, \hat{\bm{v}}_{-i}) - \min_{v_i\in\mathcal{A}}J_i^N(v_i,\hat{\bm{v}}_{-i})\right\rbrace
  = \max_{1\leq i\leq N} \left\lbrace J_i^N(\hat{v}_i, \hat{\bm{v}}_{-i}) - J_i^N(\bar{v}_i,\hat{\bm{v}}_{-i})\right\rbrace,
\end{equation}
where $\bar{v}_i$ is the best response that solves $\min_{v_i\in\mathcal{A}}J_i^N(v_i,\hat{\bm{v}}_{-i})$. We attain $\bar{v}_i$ from the following optimal control problem, while keeping other cars' strategies $\hat{\bm{v}}_{-i}$ unchanged:
\begin{align}
  &\bar{v}_i \triangleq \argmin_{v_i\in\mathcal{A}} J_i^N(v_i,\hat{\bm{v}}_{-i})=\argmin_{v_i\in\mathcal{A}}\int_0^T f(v_i(t),\rho^N_\sigma(x_i(t),t))\,dt,\label{eq:vbar1} \\
  \text{s.t. } & \dot{x}_j(t)=\hat{v}_j(t),\quad x_j(0)=x_{j,0},\quad j=1,\dots,i-1,i+1,\dots,N;\label{eq:vbar2} \\
  & \dot{x}_i(t)=v_i(t),\quad x_i(0)=x_{i,0}.\label{eq:vbar3}
\end{align}

\begin{definition}
  The \emph{accuracy} of a tuple of controls $\hat{v}_1(t),\dots,\hat{v}_N(t)$ is the $\hat{\epsilon}$ defined in Eq.~\eqref{eq:fidelity}.
\end{definition}

\subsection{Accuracy validation with numerical examples}

Summarizing Section~\ref{sec:MF}, Section~\ref{sec:sol} and Section~\ref{sec:construction}, we shall reiterate the procedure of solving an approximate equilibrium of the DG from its respective MFG in a more systematic way. The procedure of solving MFE-constructed controls and validating the accuracy of those controls is listed in Algorithm~\ref{alg:mfg2dg}. We will test some numerical examples following the procedure.

\begin{algorithm}[H]
  \centering
  \caption{Construction and Validation of MFE-constructed controls}
  \label{alg:mfg2dg}
  \begin{algorithmic}[1]
    \Require Number of cars $N$, cost function $f(u,\rho)$, space domain $[0,L]$, time horizon $[0,T]$, cars' initial positions $x_{1,0},\dots,x_{N,0}$, terminal cost $V_T(x)$, smoothing parameter $\sigma$.
    \Ensure MFE-constructed controls $\hat{v}_1(t),\dots,\hat{v}_N(t)$ and their maximal relative accuracy (MaxRA) and mean relative accuracy (MeanRA).
    \State Compute a smooth density profile $\rho_0(x)$ from cars' initial positions $x_{1,0},\dots,x_{N,0}$ by Eq.~\eqref{eq:smooth}. 
    \State Compute the MFE solution $u^*(x,t)$ on $[0,L]\times[0,T]$ from initial density $\rho_0(x)$ and terminal cost $V_T(x)$ using the algorithm proposed in Section~\ref{sec:sol}.
    \State Solve the MFE-constructed controls $\hat{v}_1(t),\dots,\hat{v}_N(t)$ from Eqs.~\eqref{eq:dyn2}\eqref{eq:dyn1}.
    \State Solve the best response strategy $\bar{v}_i(t)$ from Eqs.~\eqref{eq:vbar1}\eqref{eq:vbar2}\eqref{eq:vbar3} and obtain the improved cost $J_i^N(\bar{v}_i,\hat{\bm{v}}_{-i})$ for $i=1,\dots,N$.
    \State Compute the accuracy associated with the $i^{\text{th}}$ car $\hat{\epsilon}_i = J_i^N(\hat{v}_i, \hat{\bm{v}}_{-i}) - J_i^N(\bar{v}_i,\hat{\bm{v}}_{-i})$ for $i=1,\dots,N$.
    \State Compute the maximal relative accuracy and mean relative accuracy:
    \begin{align}
      \text{MaxRA}=\frac{\max_{i=1}^N\hat{\epsilon}_i}{\max_{i=1}^N |J_i^N(\hat{v}_i,\hat{\bm{v}}_{-i})|},\quad
      \text{MeanRA}=\frac{\frac1N\sum_{i=1}^N\hat{\epsilon}_i}{\frac1N\sum_{i=1}^N |J_i^N(\hat{v}_i,\hat{\bm{v}}_{-i})|}.
    \end{align}
  \end{algorithmic}
\end{algorithm}

In the numerical examples, we aim to convey two main messages:
\begin{enumerate}
  \item Given $N, \sigma>0$, we aim to construct an $\epsilon$-Nash equilibrium of the DG from a MFE solution and compute its accuracy $\hat{\epsilon}$.
  \item The accuracy $\hat{\epsilon}$ deceases as $N$ becomes large.
\end{enumerate}


The general set-up of the following numerical examples is similar to that used previously. We fix the length of the planning horizon $T=1$.
Each time we solve a different DG by varying the number of cars $N=21,41,61,81,101$.
For different numbers of cars, the cars' initial positions are sampled from the same initial distribution defined in Eq.~\eqref{eq:simuini} with $\rho_a=0.2$, $\rho_b=0.8$ and $\gamma=0.15L$.
For each $N$ we take the road length $L=N$ and choose the smoothing parameter $\sigma=0.05L$.
The cost functions for [MFG-NonSeparable] and [MFG-Separable] are tested.

\begin{figure}[htbp]
  \centering
  \includegraphics[width=.96\textwidth]{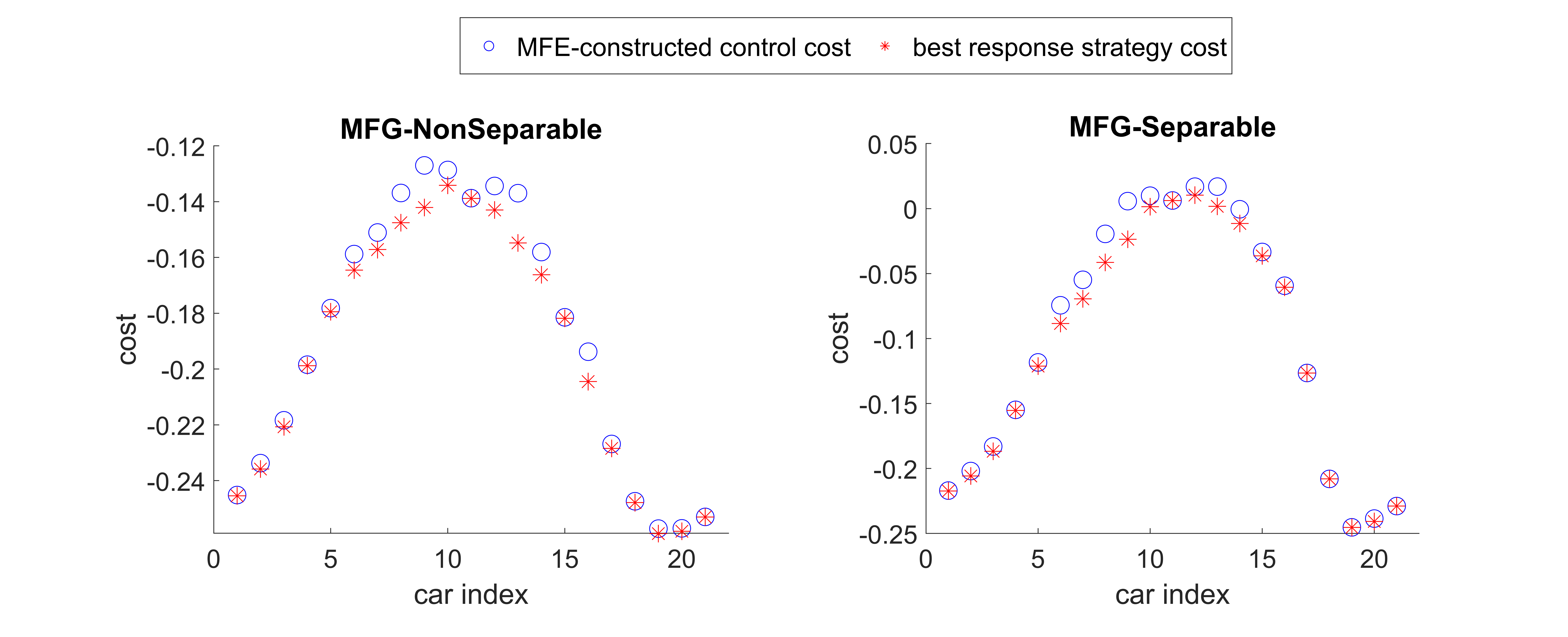}
  \caption{MFE-constructed control cost v.s. best response strategy cost, $N=21$ cars}
  \label{fig:apprvsoptim}
\end{figure}

To see the first message, Figure~\ref{fig:apprvsoptim} compares the difference between costs computed from MFE-constructed controls and those from the respective best response strategies when $N=21$. In Figure~\ref{fig:apprvsoptim}, the $x$-axis is the car's index $i$ while the $y$-axis is the cost $J_i^N$ for $i=1,2,\dots,N$. We see from the figure that for both cost functions the MFE solutions generate good approximate equilibria of DGs.

\begin{figure}[htbp]
  \centering
  \includegraphics[width=.96\textwidth]{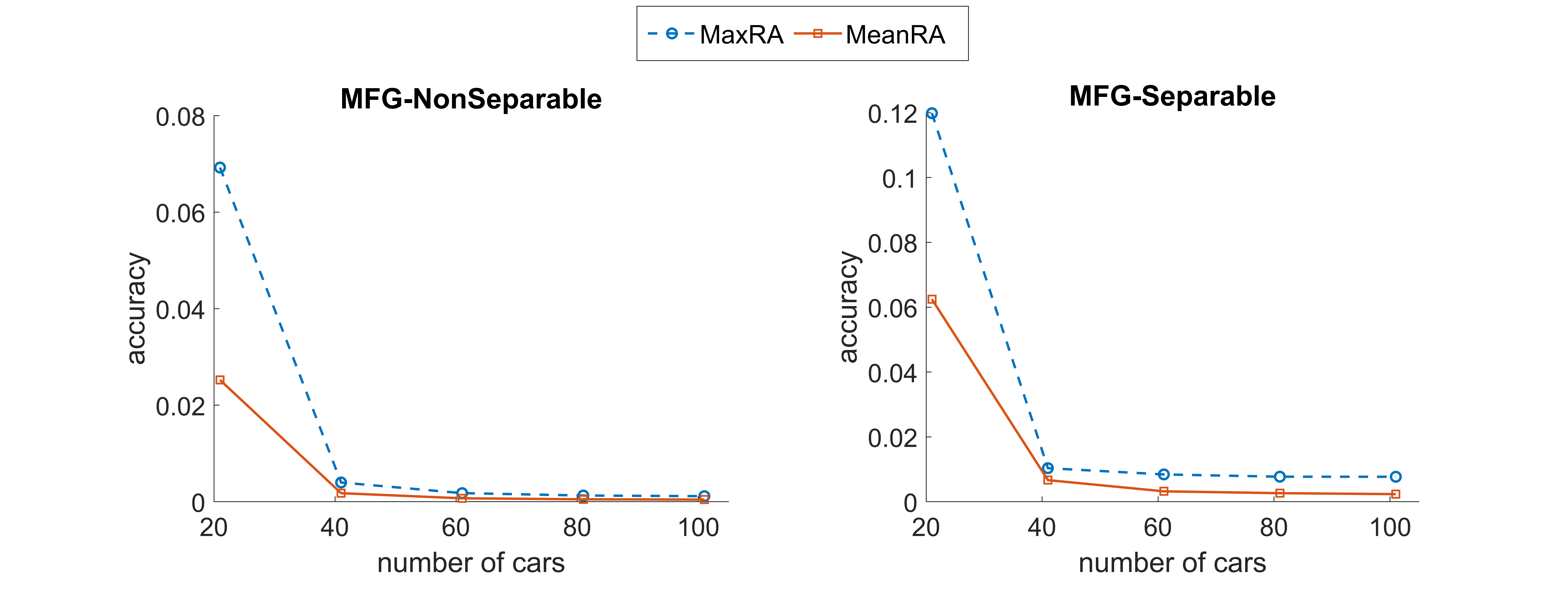}
  \caption{Accuracy v.s. Number of cars}
  \label{fig:errorvsN}
\end{figure}

To see the second message, we compute the maximal relative accuracy and mean relative accuracy of MFE-constructed controls for different numbers of cars, which is shown in Figure~\ref{fig:errorvsN}. We see from Figure~\ref{fig:errorvsN} that for both cost functions better accuracy is obtained as $N$ increases.

\subsection{Discussion}

The aforementioned procedure aims to characterize the interplay between microscopic DG and macroscopic MFG.
In particular, we can use the MFE solution to construct an $\epsilon$-Nash equilibrium of DG and show that this $\epsilon$-Nash equilibrium has desired accuracy.
It provides an efficient and scalable method to solve AVs' individual controls in a game-theoretic framework.

Note that it is challenging to find the equilibrium of DG accurately. In fact, it is not known whether or not DGE does exist and is unique.
In this paper, though, we have not discussed the solution properties of the original differential game.
There are three cases in terms of its solutions:

\begin{enumerate}
  \item DGE does not exist:
  Albeit non-existent, in practice, we still need to find a ``good enough" control for each individual AV so that every AV achieves its predefined driving objective with a reasonable performance.
  The MFE-constructed controls can be used as an approximate equilibrium of DG.
  \item DGE exists and is unique: The MFE-constructed controls provide a good initial guess for solving the accurate DGE, and the proposed method can help to characterize the upper bound of the deviation.
  \item DGE exists but is non-unique: The MFE-constructed controls can be an approximation to one DGE.
      However, the characterization of the error bound and the proposed method of finding an $\epsilon$-Nash equilibrium may become debatable.
      Such a case will be left for future research.
\end{enumerate}

\section{Conclusions and future research}\label{sec:conclude}

This paper applies the MFG to solve the continuous velocity control
problem for a system of AVs in traffic.
The MFG offers advantages over the existing individual control methods due to its scalability properties.
The proposed game-theoretic framework links micro- and macro-scale behaviors, offering insights into systematic impacts of strategic interactions among AVs from a microscopic scale.
To the best of our knowledge, this is the first study to characterize equilibrium solutions in both continuous MFGs and discrete differential games in traffic.
Unlike most of the existing studies that approximate discrete AV controls directly, we develop a game-theoretic framework from micro- to macro-scale, and then construct solutions from macro- back to micro-scale.
In particular, we first introduce the macroscopic mean field game, solve its equilibrium, construct discrete controls from the mean field equilibrium, and then validate the consistency between the constructed discrete controls and the equilibrium of the original differential game.
Our findings will help transportation engineers and planners to better predict and forecast traffic conditions when AVs reach a critical mass,
which in turn will prepare them for a smooth transition from the present to future AV-equipped transportation systems. 

This work can be generalized in several ways:
(i) As a first step to model AVs' strategic interactions under the mean field game framework, this paper only presents AVs' longitudinal velocity controls with many simplifying assumptions. The control variables, assumptions on cost functions and constraints need to be generalized.
{
For example, we will incorporate the acceleration rate into AVs' control variables}, relax the homogeneity assumptions by considering multi-class MFGs, relax assumption (A3) by incorporating the derivatives of the density into AVs' driving cost functionals and add constraints on the density and acceleration rate for the MFG.
(ii) This paper mainly focuses on deterministic, non-viscous MFGs. In the future we will discuss viscous MFGs by adding randomness in AVs' dynamics.
(iii) We rigorously show the relations between MFGs and the LWR model. However, we are only able to discuss the relationship between MFGs and traditional higher-order traffic flow models by rewriting a specific MFG system to its reduced form and find the similarity. Exploring deeper connections between MFGs and traditional higher-order models will be left for future research.
{
(iv) This paper presents three concrete cost functions and their respective MFGs to illustrate AVs' traffic flow patterns and the consistency between discrete and continuous equilibria. In the future we will explore other families of cost functions and provide more MFG examples.}

\section*{Acknowledgments}

The authors would like to thank Data Science Institute from Columbia University for providing a seed grant for this research.

\bibliographystyle{elsarticle-num}

\bibliography{survey,survey_DCL,survey-Di}

\begin{thebibliography}{100}
\expandafter\ifx\csname url\endcsname\relax
  \def\url#1{\texttt{#1}}\fi
\expandafter\ifx\csname urlprefix\endcsname\relax\def\urlprefix{URL }\fi
\expandafter\ifx\csname href\endcsname\relax
  \def\href#1#2{#2} \def\path#1{#1}\fi

\bibitem{newell1961nonlinear}
G.~F. Newell, Nonlinear effects in the dynamics of car following, Operations
  research 9~(2) (1961) 209--229.

\bibitem{gipps1981behavioural}
P.~G. Gipps, A behavioural car-following model for computer simulation,
  Transportation Research Part B: Methodological 15~(2) (1981) 105--111.

\bibitem{bando1995dynamical}
M.~Bando, K.~Hasebe, A.~Nakayama, A.~Shibata, Y.~Sugiyama, Dynamical model of
  traffic congestion and numerical simulation, Physical review E 51~(2) (1995)
  1035.

\bibitem{brackstone1999car}
M.~Brackstone, M.~McDonald, Car-following: a historical review, Transportation
  Research Part F: Traffic Psychology and Behaviour 2~(4) (1999) 181--196.

\bibitem{zhang1999mathematical}
H.~M. Zhang, A mathematical theory of traffic hysteresis, Transportation
  Research Part B: Methodological 33~(1) (1999) 1--23.

\bibitem{zhang2005car}
H.~M. Zhang, T.~Kim, A car-following theory for multiphase vehicular traffic
  flow, Transportation Research Part B: Methodological 39~(5) (2005) 385--399.

\bibitem{zheng2011applications}
Z.~Zheng, S.~Ahn, D.~Chen, J.~Laval, Applications of wavelet transform for
  analysis of freeway traffic: Bottlenecks, transient traffic, and traffic
  oscillations, Transportation Research Part B: Methodological 45~(2) (2011)
  372--384.

\bibitem{chen2012behavioral}
D.~Chen, J.~Laval, Z.~Zheng, S.~Ahn, A behavioral car-following model that
  captures traffic oscillations, Transportation research part B: methodological
  46~(6) (2012) 744--761.

\bibitem{orosz2010traffic}
G.~Orosz, R.~E. Wilson, G.~St{\'e}p{\'a}n, Traffic jams: dynamics and control
  (2010).

\bibitem{cui2017stabilizing}
S.~Cui, B.~Seibold, R.~Stern, D.~B. Work, Stabilizing traffic flow via a single
  autonomous vehicle: Possibilities and limitations, in: Intelligent Vehicles
  Symposium (IV), 2017 IEEE, IEEE, 2017, pp. 1336--1341.

\bibitem{yang2010development}
H.-H. Yang, H.~Peng, Development of an errorable car-following driver model,
  Vehicle System Dynamics 48~(6) (2010) 751--773.

\bibitem{huang2017accelerated}
Z.~Huang, H.~Lam, D.~J. LeBlanc, D.~Zhao, Accelerated evaluation of automated
  vehicles using piecewise mixture models, IEEE Transactions on Intelligent
  Transportation Systems (2017).

\bibitem{zhao2017accelerated}
D.~Zhao, H.~Lam, H.~Peng, S.~Bao, D.~J. LeBlanc, K.~Nobukawa, C.~S. Pan,
  Accelerated evaluation of automated vehicles safety in lane-change scenarios
  based on importance sampling techniques, IEEE transactions on intelligent
  transportation systems 18~(3) (2017) 595--607.

\bibitem{zhao2018accelerated}
D.~Zhao, X.~Huang, H.~Peng, H.~Lam, D.~J. LeBlanc, Accelerated evaluation of
  automated vehicles in car-following maneuvers, IEEE Transactions on
  Intelligent Transportation Systems 19~(3) (2018) 733--744.

\bibitem{liu2015safe}
C.~Liu, M.~Tomizuka, Safe exploration: Addressing various uncertainty levels in
  human robot interactions, in: 2015 American Control Conference (ACC), IEEE,
  2015, pp. 465--470.

\bibitem{liu2016enabling}
C.~Liu, M.~Tomizuka, Enabling safe freeway driving for automated vehicles, in:
  American Control Conference (ACC), 2016, IEEE, 2016, pp. 3461--3467.

\bibitem{gong2016constrained}
S.~Gong, J.~Shen, L.~Du, Constrained optimization and distributed computation
  based car following control of a connected and autonomous vehicle platoon,
  Transportation Research Part B: Methodological 94 (2016) 314--334.

\bibitem{gong2018cooperative}
S.~Gong, L.~Du, Cooperative platoon control for a mixed traffic flow including
  human drive vehicles and connected and autonomous vehicles, Transportation
  research part B: methodological 116 (2018) 25--61.

\bibitem{newell1993simplified}
G.~F. Newell, A simplified theory of kinematic waves in highway traffic, {P}art
  {I}: General theory, Transportation Research Part B: Methodological 27~(4)
  (1993) 281--287.

\bibitem{lighthill1955kinematic}
M.~J. Lighthill, G.~B. Whitham, On kinematic waves {II}. {A} theory of traffic
  flow on long crowded roads, Proc. R. Soc. Lond. A 229~(1178) (1955) 317--345.

\bibitem{zhang1998theory}
H.~M. Zhang, A theory of nonequilibrium traffic flow, Transportation Research
  Part B: Methodological 32~(7) (1998) 485--498.

\bibitem{zhang2002non}
H.~M. Zhang, A non-equilibrium traffic model devoid of gas-like behavior,
  Transportation Research Part B: Methodological 36~(3) (2002) 275--290.

\bibitem{jin2003distribution}
W.~Jin, H.~M. Zhang, On the distribution schemes for determining flows through
  a merge, Transportation Research Part B: Methodological 37~(6) (2003)
  521--540.

\bibitem{daganzo2005variational}
C.~F. Daganzo, A variational formulation of kinematic waves: basic theory and
  complex boundary conditions, Transportation Research Part B: Methodological
  39~(2) (2005) 187--196.

\bibitem{lebacque2005first}
J.-P. Lebacque, First-order macroscopic traffic flow models: Intersection
  modeling, network modeling, in: Transportation and Traffic Theory. Flow,
  Dynamics and Human Interaction. 16th International Symposium on
  Transportation and Traffic TheoryUniversity of Maryland, College Park, 2005.

\bibitem{zhou2017rolling}
Y.~Zhou, S.~Ahn, M.~Chitturi, D.~A. Noyce, Rolling horizon stochastic optimal
  control strategy for {ACC} and {CACC} under uncertainty, Transportation
  Research Part C: Emerging Technologies 83 (2017) 61--76.

\bibitem{wei2017dynamic}
Y.~Wei, C.~Avc{\i}, J.~Liu, B.~Belezamo, N.~Ayd{\i}n, P.~T. Li, X.~Zhou,
  Dynamic programming-based multi-vehicle longitudinal trajectory optimization
  with simplified car following models, Transportation research part B:
  methodological 106 (2017) 102--129.

\bibitem{li2018nonlinear}
Y.~Li, C.~Tang, S.~Peeta, Y.~Wang, Nonlinear consensus-based connected vehicle
  platoon control incorporating car-following interactions and heterogeneous
  time delays, IEEE Transactions on Intelligent Transportation Systems (2018).

\bibitem{ma2016freeway}
J.~Ma, X.~Li, S.~E. Shladover, H.~A. Rakha, X.-Y. Lu, R.~Jagannathan, D.~J.
  Dailey, Freeway speed harmonization, IEEE Trans. Intelligent Vehicles 1~(1)
  (2016) 78--89.

\bibitem{malikopoulos2018optimal}
A.~A. Malikopoulos, S.~Hong, B.~B. Park, J.~Lee, S.~Ryu, Optimal control for
  speed harmonization of automated vehicles, IEEE Transactions on Intelligent
  Transportation Systems~(99) (2018) 1--13.

\bibitem{arefizadeh2018platooning}
S.~Arefizadeh, A.~Talebpour, A platooning strategy for automated vehicles in
  the presence of speed limit fluctuations, Transportation Research Record
  2672~(20) (2018) 154--161.

\bibitem{li2018piecewise}
X.~Li, A.~Ghiasi, Z.~Xu, X.~Qu, A piecewise trajectory optimization model for
  connected automated vehicles: Exact optimization algorithm and queue
  propagation analysis, Transportation Research Part B: Methodological 118
  (2018) 429--456.

\bibitem{altan2017glidepath}
O.~D. Altan, G.~Wu, M.~J. Barth, K.~Boriboonsomsin, J.~A. Stark, Glidepath:
  Eco-friendly automated approach and departure at signalized intersections,
  IEEE Transactions on Intelligent Vehicles 2~(4) (2017) 266--277.

\bibitem{hao2018eco}
P.~Hao, G.~Wu, K.~Boriboonsomsin, M.~J. Barth, Eco-approach and departure (ead)
  application for actuated signals in real-world traffic, IEEE Transactions on
  Intelligent Transportation Systems 20~(1) (2018) 30--40.

\bibitem{yao2018trajectory}
H.~Yao, J.~Cui, X.~Li, Y.~Wang, S.~An, A trajectory smoothing method at
  signalized intersection based on individualized variable speed limits with
  location optimization, Transportation Research Part D: Transport and
  Environment 62 (2018) 456--473.

\bibitem{swaroop1996string}
D.~Swaroop, J.~K. Hedrick, String stability of interconnected systems, IEEE
  transactions on automatic control 41~(3) (1996) 349--357.

\bibitem{darbha1999intelligent}
S.~Darbha, K.~Rajagopal, Intelligent cruise control systems and traffic flow
  stability, Transportation Research Part C: Emerging Technologies 7~(6) (1999)
  329--352.

\bibitem{ioannou1993autonomous}
P.~A. Ioannou, C.-C. Chien, Autonomous intelligent cruise control, IEEE
  Transactions on Vehicular technology 42~(4) (1993) 657--672.

\bibitem{rajamani2001experimental}
R.~Rajamani, S.~E. Shladover, An experimental comparative study of autonomous
  and co-operative vehicle-follower control systems, Transportation Research
  Part C: Emerging Technologies 9~(1) (2001) 15--31.

\bibitem{van2006impact}
B.~Van~Arem, C.~J. Van~Driel, R.~Visser, The impact of cooperative adaptive
  cruise control on traffic-flow characteristics, IEEE Transactions on
  Intelligent Transportation Systems 7~(4) (2006) 429--436.

\bibitem{naus2010string}
G.~J. Naus, R.~P. Vugts, J.~Ploeg, M.~J. van~de Molengraft, M.~Steinbuch,
  String-stable {CACC} design and experimental validation: {A} frequency-domain
  approach, IEEE Transactions on vehicular technology 59~(9) (2010) 4268--4279.

\bibitem{vanderwerf2001modeling}
J.~VanderWerf, S.~Shladover, N.~Kourjanskaia, M.~Miller, H.~Krishnan, Modeling
  effects of driver control assistance systems on traffic, Transportation
  Research Record: Journal of the Transportation Research Board~(1748) (2001)
  167--174.

\bibitem{kesting2010enhanced}
A.~Kesting, M.~Treiber, D.~Helbing, Enhanced intelligent driver model to access
  the impact of driving strategies on traffic capacity, Philosophical
  Transactions of the Royal Society of London A: Mathematical, Physical and
  Engineering Sciences 368~(1928) (2010) 4585--4605.

\bibitem{schakel2010effects}
W.~J. Schakel, B.~Van~Arem, B.~D. Netten, Effects of cooperative adaptive
  cruise control on traffic flow stability, in: Intelligent Transportation
  Systems (ITSC), 2010 13th International IEEE Conference on, IEEE, 2010, pp.
  759--764.

\bibitem{ploeg2011design}
J.~Ploeg, B.~T. Scheepers, E.~Van~Nunen, N.~Van~de Wouw, H.~Nijmeijer, Design
  and experimental evaluation of cooperative adaptive cruise control, in:
  Intelligent Transportation Systems (ITSC), 2011 14th International IEEE
  Conference on, IEEE, 2011, pp. 260--265.

\bibitem{milanes2014cooperative}
V.~Milan{\'e}s, S.~E. Shladover, J.~Spring, C.~Nowakowski, H.~Kawazoe,
  M.~Nakamura, Cooperative adaptive cruise control in real traffic situations,
  IEEE Transactions on Intelligent Transportation Systems 15~(1) (2014)
  296--305.

\bibitem{milanes2014modeling}
V.~Milan{\'e}s, S.~E. Shladover, Modeling cooperative and autonomous adaptive
  cruise control dynamic responses using experimental data, Transportation
  Research Part C: Emerging Technologies 48 (2014) 285--300.

\bibitem{talebpour2016influence}
A.~Talebpour, H.~S. Mahmassani, Influence of connected and autonomous vehicles
  on traffic flow stability and throughput, Transportation Research Part C:
  Emerging Technologies 71 (2016) 143--163.

\bibitem{jin2014dynamics}
I.~G. Jin, G.~Orosz, Dynamics of connected vehicle systems with delayed
  acceleration feedback, Transportation Research Part C: Emerging Technologies
  46 (2014) 46--64.

\bibitem{levin2016multiclass}
M.~W. Levin, S.~D. Boyles, A multiclass cell transmission model for shared
  human and autonomous vehicle roads, Transportation Research Part C: Emerging
  Technologies 62 (2016) 103--116.

\bibitem{wang2015game}
M.~Wang, S.~P. Hoogendoorn, W.~Daamen, B.~van Arem, R.~Happee, Game theoretic
  approach for predictive lane-changing and car-following control,
  Transportation Research Part C: Emerging Technologies 58 (2015) 73--92.

\bibitem{Yu2018}
H.~Yu, H.~E. Tseng, R.~Langari,
  \href{https://doi.org/10.1016/j.trc.2018.01.016}{{A human-like game
  theory-based controller for automatic lane changing}}, Transportation
  Research Part C: Emerging Technologies 88~(March 2017) (2018) 140--158.
\newblock \href {https://doi.org/10.1016/j.trc.2018.01.016}
  {\path{doi:10.1016/j.trc.2018.01.016}}.
\newline\urlprefix\url{https://doi.org/10.1016/j.trc.2018.01.016}

\bibitem{RSMUP_2014__131__217_0}
R.~M. Colombo, E.~Rossi, On the micro-macro limit in traffic flow, Rendiconti
  del Seminario Matematico della Universit\`a di Padova 131 (2014) 217--236.

\bibitem{rossi2014justification}
E.~Rossi, A justification of a {LWR} model based on a follow the leader
  description., Discrete \& Continuous Dynamical Systems-Series S 7~(3) (2014).

\bibitem{di2015rigorous}
M.~Di~Francesco, M.~D. Rosini, Rigorous derivation of nonlinear scalar
  conservation laws from follow-the-leader type models via many particle limit,
  Archive for rational mechanics and analysis 217~(3) (2015) 831--871.

\bibitem{daganzo2006traffic}
C.~F. Daganzo, In traffic flow, cellular automata= kinematic waves,
  Transportation Research Part B: Methodological 40~(5) (2006) 396--403.

\bibitem{daganzo2006variational}
C.~F. Daganzo, On the variational theory of traffic flow: well-posedness,
  duality and applications (2006).

\bibitem{helbing2009derivation}
D.~Helbing, Derivation of non-local macroscopic traffic equations and
  consistent traffic pressures from microscopic car-following models, The
  European Physical Journal B 69~(4) (2009) 539--548.

\bibitem{Laval2013}
J.~A. Laval, L.~Leclercq, The {H}amilton--{J}acobi partial differential
  equation and the three representations of traffic flow, Transportation
  Research Part B: Methodological 52 (2013) 17--30.

\bibitem{jin2016equivalence}
W.-L. Jin, On the equivalence between continuum and car-following models of
  traffic flow, Transportation Research Part B: Methodological 93 (2016)
  543--559.

\bibitem{ngoduy2009continuum}
D.~Ngoduy, S.~Hoogendoorn, R.~Liu, Continuum modeling of cooperative traffic
  flow dynamics, Physica A: Statistical Mechanics and its Applications 388~(13)
  (2009) 2705--2716.

\bibitem{porfyri2015stability}
K.~Porfyri, I.~Nikolos, A.~Delis, M.~Papageorgiou, Stability analysis of a
  macroscopic traffic flow model for adaptive cruise control systems, in: ASME
  2015 International Mechanical Engineering Congress and Exposition, American
  Society of Mechanical Engineers, 2015, pp. V012T15A002--V012T15A002.

\bibitem{lighthill1952sound}
M.~J. Lighthill, On sound generated aerodynamically {I}. general theory, Proc.
  R. Soc. Lond. A 211~(1107) (1952) 564--587.

\bibitem{Payne1971}
H.~J. Payne, Model of freeway traffic and control, Mathematical Model of Public
  System (1971) 51--61.

\bibitem{Whitham2011}
G.~B. Whitham, Linear and nonlinear waves, Vol.~42, John Wiley \& Sons, 2011.

\bibitem{aw2000resurrection}
A.~Aw, M.~Rascle, Resurrection of "second order" models of traffic flow, SIAM
  journal on applied mathematics 60~(3) (2000) 916--938.

\bibitem{swaroop2001direct}
D.~Swaroop, J.~K. Hedrick, S.~B. Choi, Direct adaptive longitudinal control of
  vehicle platoons, IEEE Transactions on Vehicular Technology 50~(1) (2001)
  150--161.

\bibitem{stern2018dissipation}
R.~E. Stern, S.~Cui, M.~L. Delle~Monache, R.~Bhadani, M.~Bunting, M.~Churchill,
  N.~Hamilton, H.~Pohlmann, F.~Wu, B.~Piccoli, et~al., Dissipation of
  stop-and-go waves via control of autonomous vehicles: Field experiments,
  Transportation Research Part C: Emerging Technologies 89 (2018) 205--221.

\bibitem{swaroop1994comparision}
D.~Swaroop, J.~K. Hedrick, C.~Chien, P.~Ioannou, A comparision of spacing and
  headway control laws for automatically controlled vehicles, Vehicle system
  dynamics 23~(1) (1994) 597--625.

\bibitem{wang2016cooperative}
M.~Wang, W.~Daamen, S.~P. Hoogendoorn, B.~van Arem, Cooperative car-following
  control: Distributed algorithm and impact on moving jam features, IEEE
  Transactions on Intelligent Transportation Systems 17~(5) (2016) 1459--1471.

\bibitem{treiber2000congested}
M.~Treiber, A.~Hennecke, D.~Helbing, Congested traffic states in empirical
  observations and microscopic simulations, Physical review E 62~(2) (2000)
  1805.

\bibitem{kesting2008adaptive}
A.~Kesting, M.~Treiber, M.~Sch{\"o}nhof, D.~Helbing, Adaptive cruise control
  design for active congestion avoidance, Transportation Research Part C:
  Emerging Technologies 16~(6) (2008) 668--683.

\bibitem{wu2017emergent}
C.~Wu, A.~Kreidieh, E.~Vinitsky, A.~M. Bayen, Emergent behaviors in
  mixed-autonomy traffic, in: Conference on Robot Learning, 2017, pp. 398--407.

\bibitem{talebpour2018effect}
A.~Talebpour, H.~S. Mahmassani, S.~H. Hamdar, et~al., Effect of information
  availability on stability of traffic flow: Percolation theory approach,
  Transportation Research Part B: Methodological 117~(PB) (2018) 624--638.

\bibitem{wu2018stabilizing}
C.~Wu, A.~M. Bayen, A.~Mehta, Stabilizing traffic with autonomous vehicles, in:
  2018 IEEE International Conference on Robotics and Automation (ICRA), IEEE,
  2018, pp. 1--7.

\bibitem{qin2013digital}
W.~B. Qin, G.~Orosz, Digital effects and delays in connected vehicles: linear
  stability and simulations, in: ASME 2013 Dynamic Systems and Control
  Conference, American Society of Mechanical Engineers, 2013, pp.
  V002T30A001--V002T30A001.

\bibitem{qin2017scalable}
W.~B. Qin, G.~Orosz, Scalable stability analysis on large connected vehicle
  systems subject to stochastic communication delays, Transportation Research
  Part C: Emerging Technologies 83 (2017) 39--60.

\bibitem{jin2018connected}
I.~G. Jin, G.~Orosz, Connected cruise control among human-driven vehicles:
  Experiment-based parameter estimation and optimal control design,
  Transportation Research Part C: Emerging Technologies 95 (2018) 445--459.

\bibitem{jin2018experimental}
I.~G. Jin, S.~S. Avedisov, C.~R. He, W.~B. Qin, M.~Sadeghpour, G.~Orosz,
  Experimental validation of connected automated vehicle design among
  human-driven vehicles, Transportation research part C: emerging technologies
  91 (2018) 335--352.

\bibitem{wang2014rolling2}
M.~Wang, W.~Daamen, S.~P. Hoogendoorn, B.~van Arem, Rolling horizon control
  framework for driver assistance systems. {P}art {II}: Cooperative sensing and
  cooperative control, Transportation research part C: emerging technologies 40
  (2014) 290--311.

\bibitem{wang2014rolling1}
M.~Wang, W.~Daamen, S.~P. Hoogendoorn, B.~van Arem, Rolling horizon control
  framework for driver assistance systems. {P}art {I}: Mathematical formulation
  and non-cooperative systems, Transportation research part C: emerging
  technologies 40 (2014) 271--289.

\bibitem{dreves2018generalized}
A.~Dreves, M.~Gerdts, A generalized {N}ash equilibrium approach for optimal
  control problems of autonomous cars, Optimal Control Applications and Methods
  39~(1) (2018) 326--342.

\bibitem{li2018game}
N.~Li, D.~W. Oyler, M.~Zhang, Y.~Yildiz, I.~Kolmanovsky, A.~R. Girard, Game
  theoretic modeling of driver and vehicle interactions for verification and
  validation of autonomous vehicle control systems, IEEE Transactions on
  control systems technology 26~(5) (2018) 1782--1797.

\bibitem{sadigh2016planning}
D.~Sadigh, S.~Sastry, S.~A. Seshia, A.~D. Dragan, Planning for autonomous cars
  that leverage effects on human actions., in: Robotics: Science and Systems,
  2016.

\bibitem{lazar2018maximizing}
D.~A. Lazar, K.~Chandrasekher, R.~Pedarsani, D.~Sadigh, Maximizing road
  capacity using cars that influence people, arXiv preprint arXiv:1807.04414
  (2018).

\bibitem{Talebpour2015}
A.~Talebpour, H.~S. Mahmassani, S.~H. Hamdar, Modeling lane-changing behavior
  in a connected environment: {A} game theory approach, Transportation Research
  Procedia 7 (2015) 420--440.
\newblock \href {https://doi.org/10.1016/j.trpro.2015.06.022}
  {\path{doi:10.1016/j.trpro.2015.06.022}}.

\bibitem{yoo2012stackelberg}
J.~H. Yoo, R.~Langari, Stackelberg game based model of highway driving, in:
  ASME 2012 5th Annual Dynamic Systems and Control Conference joint with the
  JSME 2012 11th Motion and Vibration Conference, American Society of
  Mechanical Engineers, 2012, pp. 499--508.

\bibitem{yoo2013stackelberg}
J.~H. Yoo, R.~Langari, A {S}tackelberg game theoretic driver model for merging,
  in: ASME 2013 Dynamic Systems and Control Conference, American Society of
  Mechanical Engineers, 2013, pp. V002T30A003--V002T30A003.

\bibitem{kim2014game}
C.~Kim, R.~Langari, Game theory based autonomous vehicles operation,
  International Journal of Vehicle Design 65~(4) (2014) 360--383.

\bibitem{yu2018human}
H.~Yu, H.~E. Tseng, R.~Langari, A human-like game theory-based controller for
  automatic lane changing, Transportation Research Part C: Emerging
  Technologies 88 (2018) 140--158.

\bibitem{wu2017flow}
C.~Wu, A.~Kreidieh, K.~Parvate, E.~Vinitsky, A.~M. Bayen, Flow: Architecture
  and benchmarking for reinforcement learning in traffic control, arXiv
  preprint arXiv:1710.05465 (2017).

\bibitem{wu2017framework}
C.~Wu, K.~Parvate, N.~Kheterpal, L.~Dickstein, A.~Mehta, E.~Vinitsky, A.~M.
  Bayen, Framework for control and deep reinforcement learning in traffic, in:
  Intelligent Transportation Systems (ITSC), 2017 IEEE 20th International
  Conference on, IEEE, 2017, pp. 1--8.

\bibitem{lasry2007mean}
J.-M. Lasry, P.-L. Lions, Mean field games, Japanese Journal of Mathematics
  2~(1) (2007) 229--260.

\bibitem{huang2006large}
M.~Huang, R.~P. Malham{\'e}, P.~E. Caines, et~al., Large population stochastic
  dynamic games: closed-loop {M}ckean-{V}lasov systems and the {N}ash certainty
  equivalence principle, Communications in Information \& Systems 6~(3) (2006)
  221--252.

\bibitem{gueant2011mean}
O.~Gu{\'e}ant, J.-M. Lasry, P.-L. Lions, Mean field games and applications, in:
  Paris-Princeton lectures on mathematical finance 2010, Springer, 2011, pp.
  205--266.

\bibitem{lachapelle2010computation}
A.~Lachapelle, J.~Salomon, G.~Turinici, Computation of mean field equilibria in
  economics, Mathematical Models and Methods in Applied Sciences 20~(04) (2010)
  567--588.

\bibitem{djehiche2016mean}
B.~Djehiche, A.~Tcheukam, H.~Tembine, Mean-field-type games in engineering,
  arXiv preprint arXiv:1605.03281 (2016).

\bibitem{couillet2012electrical}
R.~Couillet, S.~M. Perlaza, H.~Tembine, M.~Debbah, Electrical vehicles in the
  smart grid: A mean field game analysis, IEEE Journal on Selected Areas in
  Communications 30~(6) (2012) 1086--1096.

\bibitem{degond2014large}
P.~Degond, J.-G. Liu, C.~Ringhofer, Large-scale dynamics of mean-field games
  driven by local nash equilibria, Journal of Nonlinear Science 24~(1) (2014)
  93--115.

\bibitem{lachapelle2011mean}
A.~Lachapelle, M.-T. Wolfram, On a mean field game approach modeling congestion
  and aversion in pedestrian crowds, Transportation research part B:
  methodological 45~(10) (2011) 1572--1589.

\bibitem{burger2013mean}
M.~Burger, M.~Di~Francesco, P.~Markowich, M.-T. Wolfram, Mean field games with
  nonlinear mobilities in pedestrian dynamics, arXiv preprint arXiv:1304.5201
  (2013).

\bibitem{kachroo2016inverse}
P.~Kachroo, S.~Agarwal, S.~Sastry, Inverse problem for non-viscous mean field
  control: Example from traffic, IEEE Transactions on Automatic Control 61~(11)
  (2016) 3412--3421.

\bibitem{chevalier2015micro}
G.~Chevalier, J.~Le~Ny, R.~Malham{\'e}, A micro-macro traffic model based on
  mean-field games, in: 2015 American Control Conference (ACC), IEEE, 2015, pp.
  1983--1988.

\bibitem{huang2019stabilizing}
K.~Huang, X.~Di, Q.~Du, X.~Chen, Stabilizing traffic via autonomous vehicles: A
  continuum mean field game approach, in: 2019 IEEE Intelligent Transportation
  Systems Conference (ITSC), IEEE, 2019, pp. 3269--3274.

\bibitem{huang2020scalable}
K.~Huang, X.~Di, Q.~Du, X.~Chen, Scalable traffic stability analysis in
  mixed-autonomy using continuum models, Transportation research part C:
  emerging technologies 111 (2020) 616--630.

\bibitem{benamou2017variational}
J.-D. Benamou, G.~Carlier, F.~Santambrogio, Variational mean field games, in:
  Active Particles, Volume 1, Springer, 2017, pp. 141--171.

\bibitem{basar1999dynamic}
T.~Basar, G.~J. Olsder, Dynamic noncooperative game theory, Vol.~23, Siam,
  1999.

\bibitem{cardaliaguet2010notes}
P.~Cardaliaguet, Notes on mean field games, Tech. rep. (2010).

\bibitem{MFGfish}
L.~N. Hoang, The new big fish called mean-field game theory,
  \url{http://www.science4all.org/article/mean-field-games/}, [Online; accessed
  1.28.2020] (2014).

\bibitem{bardi2008optimal}
M.~Bardi, I.~Capuzzo-Dolcetta, Optimal control and viscosity solutions of
  Hamilton-Jacobi-Bellman equations, Springer Science \& Business Media, 2008.

\bibitem{shephard1974law}
R.~W. Shephard, R.~F{\"a}re, The law of diminishing returns, in: Production
  theory, Springer, 1974, pp. 287--318.

\bibitem{imprialou2016re}
M.-I.~M. Imprialou, M.~Quddus, D.~E. Pitfield, D.~Lord, Re-visiting
  crash--speed relationships: A new perspective in crash modelling, Accident
  Analysis \& Prevention 86 (2016) 173--185.

\bibitem{cardaliaguet2015weak}
P.~Cardaliaguet, Weak solutions for first order mean field games with local
  coupling, in: Analysis and geometry in control theory and its applications,
  Springer, 2015, pp. 111--158.

\bibitem{richards1956shock}
P.~I. Richards, Shock waves on the highway, Operations research 4~(1) (1956)
  42--51.

\bibitem{ambrose2018existence}
D.~M. Ambrose, Existence theory for non-separable mean field games in sobolev
  spaces, arXiv preprint arXiv:1807.02223 (2018).

\bibitem{rockafellar2015convex}
R.~T. Rockafellar, Convex analysis, Princeton university press, 2015.

\bibitem{benamou2015augmented}
J.-D. Benamou, G.~Carlier, Augmented {L}agrangian methods for transport
  optimization, mean field games and degenerate elliptic equations, Journal of
  Optimization Theory and Applications 167~(1) (2015) 1--26.

\bibitem{chow2018algorithm}
Y.~T. Chow, W.~Li, S.~Osher, W.~Yin, Algorithm for {H}amilton-{J}acobi
  equations in density space via a generalized {H}opf formula, arXiv preprint
  arXiv:1805.01636 (2018).

\bibitem{achdou2010mean}
Y.~Achdou, I.~Capuzzo-Dolcetta, Mean field games: numerical methods, SIAM
  Journal on Numerical Analysis 48~(3) (2010) 1136--1162.

\bibitem{achdou2012mean}
Y.~Achdou, F.~Camilli, I.~Capuzzo-Dolcetta, Mean field games: numerical methods
  for the planning problem, SIAM Journal on Control and Optimization 50~(1)
  (2012) 77--109.

\bibitem{achdou2012iterative}
Y.~Achdou, V.~Perez, Iterative strategies for solving linearized discrete mean
  field games systems., Networks \& Heterogeneous Media 7~(2) (2012).

\bibitem{leveque2002finite}
R.~J. LeVeque, Finite volume methods for hyperbolic problems, Vol.~31,
  Cambridge university press, 2002.

\bibitem{hackbusch2013multi}
W.~Hackbusch, Multi-grid methods and applications, Vol.~4, Springer Science \&
  Business Media, 2013.

\bibitem{golub2012matrix}
G.~H. Golub, C.~F. Van~Loan, Matrix computations, Vol.~3, JHU press, 2012.

\bibitem{saad1986gmres}
Y.~Saad, M.~H. Schultz, {GMRES: A generalized minimal residual algorithm for
  solving nonsymmetric linear systems}, SIAM Journal on scientific and
  statistical computing 7~(3) (1986) 856--869.

\bibitem{siebel2006fundamental}
F.~Siebel, W.~Mauser, On the fundamental diagram of traffic flow, SIAM Journal
  on Applied Mathematics 66~(4) (2006) 1150--1162.

\bibitem{seibold2012constructing}
B.~Seibold, M.~R. Flynn, A.~R. Kasimov, R.~R. Rosales, Constructing set-valued
  fundamental diagrams from jamiton solutions in second order traffic models,
  arXiv preprint arXiv:1204.5510 (2012).

\bibitem{strogatz2018nonlinear}
S.~H. Strogatz, Nonlinear dynamics and chaos: with applications to physics,
  biology, chemistry, and engineering, CRC Press, 2018.

\end{thebibliography}
\end{document}